\documentclass{siamart0516}
\pdfoutput=1

\usepackage{lipsum}
\usepackage{amsfonts}
\usepackage{graphicx}
\usepackage{epstopdf}
\usepackage{algorithmic}
\ifpdf
  \DeclareGraphicsExtensions{.eps,.pdf,.png,.jpg}
\else
  \DeclareGraphicsExtensions{.eps}
\fi

\usepackage{bm}
\usepackage{color}
\usepackage{url}
\usepackage{hyperref}
\usepackage{enumerate}
\usepackage{bbm}
\newcommand{\ud}{\mbox{d}}

\newcommand{\TheTitle}{Fractional Euler limits and their applications} 
\newcommand{\TheAuthors}{S. MacNamara, B. Henry, and  W. McLean}

\headers{\TheTitle}{\TheAuthors}

\title{{\TheTitle}}

\author{
  Shev MacNamara
  \and
  Bruce Henry
  \and
  William McLean\thanks{The School of Mathematics and Statistics, University of New South Wales (UNSW) (\email{s.macnamara@unsw.edu.au}, \email{b.henry@unsw.edu.au}, \email{w.mclean@unsw.edu.au})}}

\usepackage{amsopn}

\ifpdf
\hypersetup{
  pdftitle={\TheTitle},
  pdfauthor={\TheAuthors}
}
\fi

\begin{document}

\maketitle

\begin{abstract}
Generalisations of the classical Euler formula to the setting of fractional calculus are discussed.
Compound interest and fractional compound interest serve as motivation.
Connections to fractional master equations are highlighted.
An application to the Schl\"ogl reactions with Mittag-Leffler waiting times is described.
\end{abstract}



\section{\label{sec:Introduction}Introduction}

Euler's famous limit formula states, as $n \rightarrow \infty$,
\begin{eqnarray}
  \left(1+ \frac{1}{n} \right)^n & \longrightarrow & \exp(1)     = e \approx 2.7183 \ldots
  \label{eq:Euler:limit}
\end{eqnarray}
Euler  (\textit{Introductio}, 1748) was motivated by quandaries such as ``if a man borrowed 400,000 florins at the usurious rate of five percent annual interest \ldots'' \cite{AnalysisByItsHistoryHairerWanner}.
Indeed, the special number $e$ of calculus and possibly also the formula \eqref{eq:Euler:limit} may have been discovered in this way in 1683, although not by Euler; it was another Swiss mathematician, Bernoulli, in connection to his own studies of compound interest.
Others certainly made contributions, including Mercator's earlier work on the logarithm.
Euler's limit is one way to introduce exponential growth: as the discrete compounding intervals become more frequent they tend to a continuous exponential limit.
The wider principle is that \textit{a discrete process limits to a continuum object.}

Such examples of compounding processes are \textit{memoryless} in the sense that the future is conditionally independent of the past, given the present.
Other processes have memory and depend on their history.
The \textit{fractional calculus} \cite{Podlubny1999,BerkowitzEderyScher2013,Henry2006,Angstmann2015,KenkreMontrollShlesinger73,MagdziarzWernon2007,MetzlerKlafter2000,McLean2010,McLeanThomee2004,MontrollWeiss65,Fedotov:2010aa,Kozubowski2000,BelinskiyKozubowski2000,Klafter:1987aa}, unknown to Euler (though his ourve is related), offers a mathematical framework for such processes.
That involves a generalisation of the derivative to a fractional derivative, alongside which the continuous exponential function is generalised to a continuous \textit{Mittag-Leffler function}\footnote{$E_{\alpha}(z)$ in \eqref{eq:mittag:leffler} is the one parameter Mittag-Leffler function, whose namesake completed his thesis at Uppsala University, while the two-parameter Mittag-Leffler function, $E_{\alpha,\beta}(z)$, which we touch on in \eqref{eq:waiting:time:mittag:leffler}, was introduced by Anders Wiman who also shares a connection to Uppsala \cite{MittagLefflerFunctionReviewArticle,MittagLeffler1903Original1,MittagLeffler1903Original2}.}:
\begin{equation}
 E_{\alpha}(z) = \sum_{k=0}^{\infty} \frac{z^{k}}{\Gamma(\alpha k +1)} .
 \label{eq:mittag:leffler}
\end{equation}
Here $\alpha $ smoothly interpolates between the usual calculus at one extreme ($\alpha=1$), and  an ever `more fractional calculus' towards the other extreme ($\alpha \rightarrow 0$).
In this article, always $0 \le \alpha \le 1$.
The familiar power series for the exponential, $ E_{1}(t)=1+ t/1!+ t^2/2! + \ldots $, is recovered when $\alpha=1$.
Fractional processes require specialised numerical algorithms \cite{McLeanSIAMNumericalAnalysis2013,YangMoroneyBurrageTurnerLiu2011} and we use codes of Garrappa \cite{Garrappa2015} for \eqref{eq:mittag:leffler}.
Table~\ref{tab:exp:mittag:leffler} collects various representations of these functions.
A missing entry suggests a question.
\textit{What might be the fractional generalisation of the discrete limiting process in Euler's famous limit formula?}
This question will eventually lead us to the \textit{resolvent} but we begin with more elementary approaches.

\begin{table}[ht]
\caption{Comparing representations of exponential functions and of Mittag-Leffler functions.}
\begin{center}
\begin{tabular}{|c|c|c|}
  \hline  & & \\
& Exponential \;\;  $e^{t}$ \;\; & Mittag-Leffler \;\; $E_{\alpha}(t) $ \;\;\; \\
\hline & & \\
Taylor Series & \( \displaystyle   \sum_{n=0}^{\infty} \dfrac{t^n}{n !} \) &  \( \displaystyle   \sum_{k=0}^{\infty} \dfrac{t^{k}}{\Gamma(\alpha k +1)} \) \\
& & \\
Cauchy  Integral & \( \displaystyle   \frac{1}{2 \pi i}  \int_{\mathcal{C}} e^{z}   \frac{1}{z -t}   \ud z  \) &  \( \displaystyle  \frac{1}{2 \pi i}  \int_{\mathcal{C}} e^{z} \frac{z^{\alpha-1}}{z^{\alpha}-t} \ud z. \) \\
& & \\
\textbf{Euler Limit} & \( \displaystyle \lim_{n \rightarrow \infty} \left(1+\dfrac{t}{n}\right)^n \) & {\large \textbf{?} }\\
\hline
\end{tabular}
\end{center}
\label{tab:exp:mittag:leffler}
\end{table}%

\subsection{\label{sec:finite:difference:euler:limit}A recursive formulation}
Revisit Euler's  formula, recursively\footnote{Various interpretations of compound interest are possible.
A different candidate, in which the interest rate is not fixed, could come by setting $r_j = t/j$ and replacing \eqref{eq:recursive:Euler:limit} with $y_j = y_{j-1} (1+r_j)$.
  A solution is possible in terms of special functions, $y_n = y_0  \Gamma(n+1+t) (\Gamma(n+1) \Gamma(1+t))^{-1} = (t B(t,n+1))^{-1}$ (where $\partial / \partial x B(x,y) = B(x,y) (\Psi(x) - \Psi(x+y))$ and $\Psi(x) = \ud / \ud x \ln \Gamma(x) = \Gamma ' (x)/\Gamma(x)$ is the digamma function), although these alternatives are not further explored here.}
 \begin{eqnarray}
  y_0 &=& 1 , \nonumber \\
  y_{j} &=& \left(1 + h \right) y_{j-1}.
 \label{eq:recursive:Euler:limit}
\end{eqnarray}
It is common to allow a time parameter $t$ in Euler's limit formula \eqref{eq:Euler:limit}, which often appears as
$
  (1+t/n)^n \rightarrow \exp(t) \; \textrm{as  }  n \rightarrow \infty.
$
This is accommodated in \eqref{eq:recursive:Euler:limit} by setting the step size to be
$
h=t/n ,
$
for fixed $t$.
The limit $n \rightarrow \infty$ is the same as  $h \rightarrow 0$, with $t$ held constant.
This recursive formulation certainly has the  property that we expect:
$
y_n \rightarrow \exp(t)  \, \textrm{as  }  n \rightarrow \infty.
$

In fact, \eqref{eq:recursive:Euler:limit} is precisely the Euler method for approximating the simple differential equation $dy/dt= y$ with familiar solution $y(t)=\exp(t)$ when $y(0)=1$.
Finite differences of the continuous equation $dy/dt = y$ lead to the discrete approximation
\begin{equation}
\frac{y_{j}-y_{j-1}}{h} = y_{j-1}.
 \label{eq:Euler:finite:difference}
\end{equation}
Rearranging yields the Euler formula.
Comparing Taylor series shows that the local error over a single time step of size $h$ is $\mathcal{O}(h^2)$, but errors accumulate over the many steps it takes to reach $y_n$ so by the last step the global error, $y_n - \exp(t)$, is $\mathcal{O}(h)$.
This is first order accuracy: convergence is slow and Euler's formula is not usually a good numerical choice.
Nevertheless, it remains pertinent to compound interest.

\textit{Stability}\footnote{Our one-step method is too simple but in more general settings the wider significance of this property arises in a fundamental \textit{meta}-theorem of numerical analysis: \textit{stability and consistency imply convergence}. 
This may be attributed to Lax in a setting of linear PDEs while for our ODE setting it may be attributed to Dahlquist, who also shares an Uppsala connection and was an academic grandson of Anders Wiman  \cite{LaxRic56,AriehIserlesBook2009}.} is another important property of a numerical method that is related to compound interest.
For example, for the stable continuous equation $dy / dt = -y$ with solution in terms of exponential decay $y(t)=e^{-t} y(0)$, stability requirements place a restriction on the step-size, $|1-h|<1$, of the \textit{explicit} Euler \textit{forward difference} construction $y_j = (1-h)y_{j-1} $.
Stability in relation to backward differences arises later in \eqref{eq:GL:fractional:Euler:limit} and \eqref{eq:euler:limit:forward:backward}.
Errors accumulate in a way that is analogous to compound interest on a bank account.
Errors made early contribute more to the final error.
Ensuring numerical stability is tantamount to ensuring that the `compound interest' on those local errors does not grow too usuriously.
\begin{figure}[ht!]
\centering
\begin{tabular}{cc}
\includegraphics[scale=0.3,trim=15mm 70mm 25mm 60mm, clip]{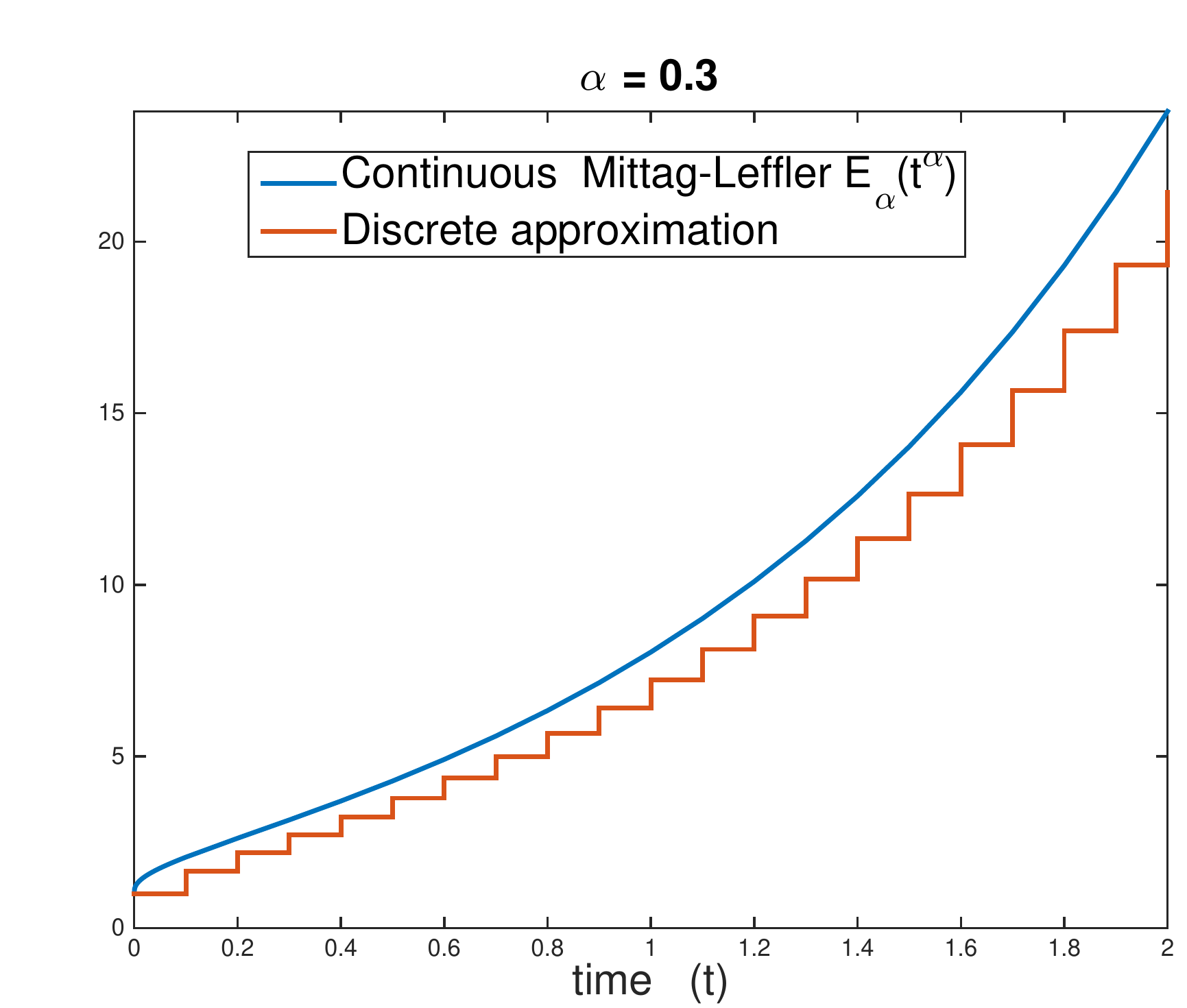} &
\includegraphics[scale=0.3,trim=15mm 70mm 25mm 70mm, clip]{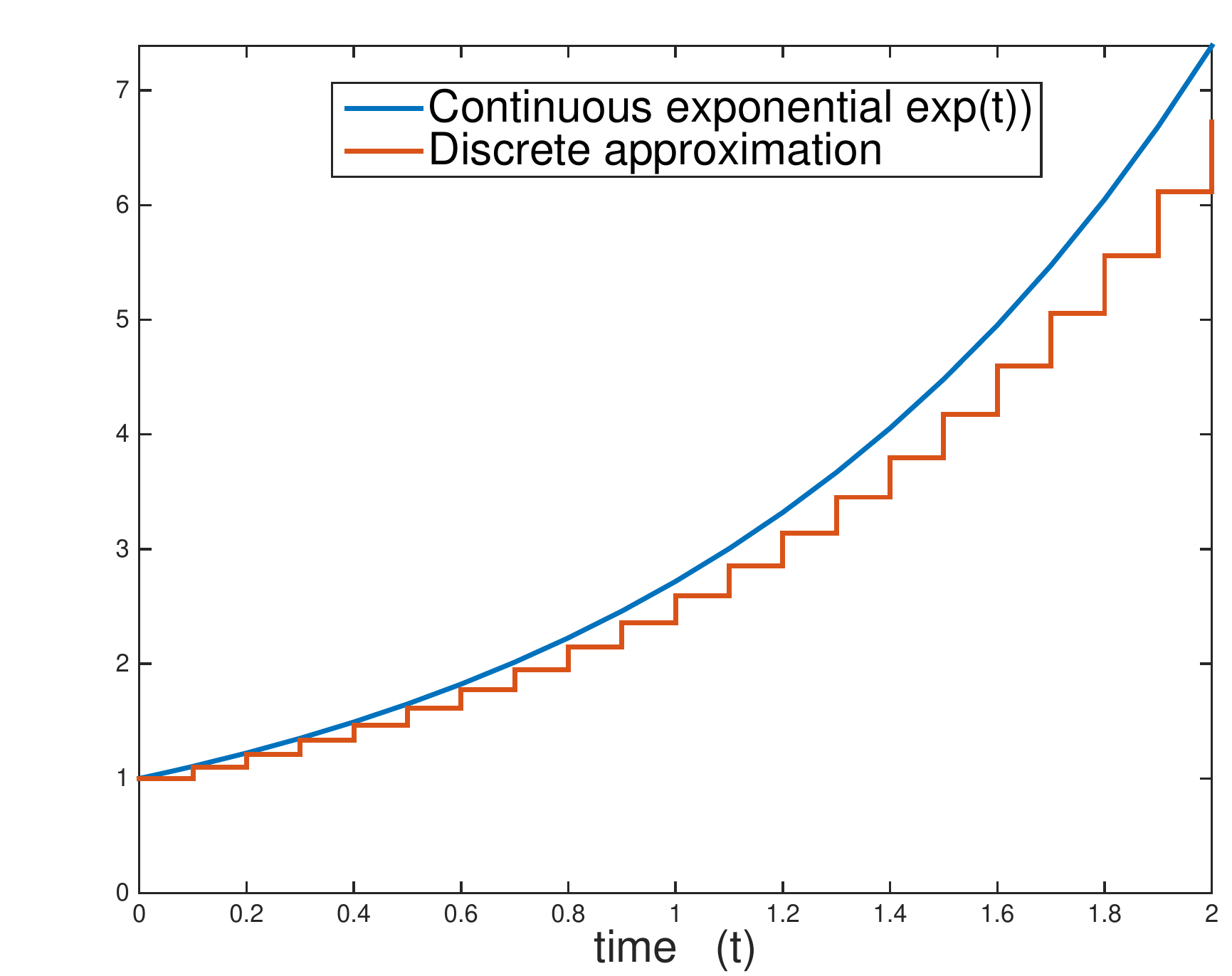}
\end{tabular}
\caption{\textit{Left}: Fractional compound interest is a Mittag-Leffler function when interest is compounded continuously, or it is the fractional generalisation proposed here in \eqref{eq:fractional:Euler:limit} when interest is compounded in discrete time steps.
 \textit{Right}: For reference, also shown is the usual interpretation of compound interest as an exponential function when compounded continuously, or as the Euler formula \eqref{eq:Euler:limit}  when compounded discretely.
 Notice the difference in scales.}
\label{fig:numerical:solution}
\end{figure}

\subsection{\label{sec:finite:difference:fractional:euler:limit}A candidate for a fractional Euler formula}
A first candidate answer to our motivating question from the Introduction could be, again recursively,
 \begin{eqnarray}
  y_0 &=& 1 ,  \nonumber  \\
 y_{j} &=& \left(1 + h^\alpha \Gamma(1 - \alpha)\, \right) y_{j-1} + \frac{y_{j-2} - y_{j-1}}{2^\alpha} + \dots \, + \frac{y_{0} - y_{1}}{j^\alpha} .
 \label{eq:fractional:Euler:limit}
\end{eqnarray}
Here $\Gamma(z)= \int_0^\infty x^{z-1} \exp(-x) \ud x$ is the Gamma function that generalises the usual factorial function $n!$ from the integers to complex numbers.
To arrive at \eqref{eq:fractional:Euler:limit}, generalise the simple finite difference construction \eqref{eq:Euler:finite:difference}  that led to Euler's limit.

Begin by defining the Caputo fractional derivative $D_t^{\alpha}$  of order $\alpha$ via its action on a smooth function
\begin{equation}
D_t^{\alpha} f(t) \equiv \frac{1}{\Gamma(1-\alpha)} \int_0^{t} \frac{f'(s)}{(t-s)^{\alpha}} \ud s
\label{eq:Caputo:derivative:definition}
\end{equation}
and introduce a time-fractional generalisation of our simple differential equation
 \begin{equation}
D_t^{\alpha} y = \lambda y \qquad \textrm{\;\;\;\;\;\;\;\;\;\; with solution \;\;\; } \quad y(t) = E_{\alpha}(\lambda t^{\alpha}) y(0).
 \label{eq:fractional:ode}
 \end{equation}
This fractional analogue offers a continuous target for a discrete process.
Set $\lambda=1$ for simplicity and difference both sides of $D_t^{\alpha} y = \lambda y$.
By quadrature on the integral in the Caputo derivative, this procedure results in
 \begin{equation}
\frac{1}{\Gamma(1-\alpha)} h \left( \frac{(y_1-y_0)/h}{(jh)^{\alpha}} + \ldots + \frac{(y_{j}-y_{j-1})/h}{(h)^{\alpha}}  \right)  = y_{j-1} . \nonumber
 \label{eq:fractional:discrete:recursion}
 \end{equation}
Rearranging yields the proposal \eqref{eq:fractional:Euler:limit}. 
Importantly,  \eqref{eq:fractional:Euler:limit} converges to the solution of the fractional differential equation in terms of the Mittag-Leffler function.
For $0< \alpha <1$,
$y_n \rightarrow E_{\alpha}(t^{\alpha})$,
 as expected.

  Figure~\ref{fig:numerical:solution} compares the usual notion of compound interest with what might be named `fractional compound interest.'
Concerning this candidate \eqref{eq:fractional:Euler:limit}, we note:
 $(i)$ The powers of $(1+h)$ present in the Euler limit are generalised to powers of $\left(1 + h^\alpha \Gamma(1 \!-
 \!\! \alpha)\! \right)$ in \eqref{eq:fractional:Euler:limit};
 $(ii)$ Convergence of the Euler limit is slow, but convergence of the fractional generalisation is even slower.
  Convergence also depends on $\alpha$.
The singularity in the gamma function $\Gamma(1-\alpha)$ as $\alpha \rightarrow 1$ is one numerical issue --- this recursion is not a good numerical method for computing the Mittag-Leffler function;
$(iii)$ The Euler limit  \eqref{eq:recursive:Euler:limit} is memoryless, requiring only the present value $y_{j-1}$ to advance to the next value of $y_{j}$, but the fractional generalisation \eqref{eq:fractional:Euler:limit} requires the whole history of values $y_0, \dots, y_{j-1}$ in order to advance;
$(iv)$ Undergraduate calculus textbooks proceed via logarithms to derive the Euler limit, so a different approach to generalising the Euler limit, not explored here, could come via a fractional generalisation of the logarithm
$ \log(t)=\int_1^t \frac{1}{u} \, \ud u$.

\textit{Fractional decay.}
A fractional generalisation of the usual exponential decay process is modelled by
$
D_t^{\alpha} y = \bm{-} y(t),
$
with solution $y(t) = E_{\alpha}(\bm{-} t^{\alpha}) y(0)$.
The same approach to discretization that led to \eqref{eq:fractional:Euler:limit} now leads to  the recursive formulation: $y_0 = 1$, and
\begin{equation}
 y_{j}  = \left(1 \bm{-} h^\alpha \Gamma(1 - \alpha)\, \right) y_{j-1} + \frac{y_{j-2} - y_{j-1}}{2^\alpha} + \dots \, + \frac{y_{0} - y_{1}}{j^\alpha}.
 \label{eq:fractional:Euler:limit:minus}
 \end{equation}
This is the counterpart to \eqref{eq:fractional:Euler:limit} when the argument to the Mittag-Leffler function is negative.
Apart from a minus sign in the first term in parentheses, it is otherwise identical with \eqref{eq:fractional:Euler:limit}: the `memory terms'  (connected to a memory function later in \eqref{eq:convolution:memory:function}) have the same form in both growth and decay versions.

\subsection{\label{sec:Grunwald:Letnikov} A Gr\"unwald-Letnikov approach}
The binomial  term $(1+h)^n$ that appears in the Euler limit suggests another approach to its generalisation could come via the Gr\"unwald-Letnikov (GL) derivative, which defines fractional derivatives via a fractional generalisation of a binomial-like expansion.
With this  in mind,  another way to express the same fractional  model $D_t^{\alpha} y = \lambda y $ in \eqref{eq:fractional:ode} is
 \begin{equation}
y(t) = y(0) + \lambda I^{\alpha} y,
 \label{eq:GL:fractional:integral}
 \end{equation}
where $I^{\alpha}$ is the GL fractional integral operator\footnote{There are different conventions for defining fractional \textit{derivatives}, which vary in details of how they handle initial conditions, but there is little distinction in relation to fractional integrals.}.
(Indeed, one way to formalise the notion of solution to $D_t^{\alpha} y = \lambda y $ is as a solution of the Volterra integral equation $y(t) = y(0) +  \int_0^t k(t-u) \lambda y(u) \ud u$ where $k(t) = t^{\alpha-1}/ \Gamma(\alpha)$ is a memory-like function.
This integral displays the causal nature of the problem: the solution at the present time can be computed using only information from the past, without needing information from the future.)
Hence \eqref{eq:GL:fractional:integral} is another representation of the Mittag-Leffler function.

A Gr\"unwald-Letnikov approach to fractional calculus is attractive because the construction is discrete from the beginning (analogous to the way that the usual calculus begins, as continuous limits of finite-difference constructions), unlike other approaches (such as  Riemann-Louiville), which begin with continuous integer-order integrals or derivatives.
Set $h=t/n$, and replace the continuous integral $I^\alpha$ in  \eqref{eq:GL:fractional:integral}  by the discrete  GL construction that led to it to obtain
$
y_n = y_0 + h^{\alpha} \sum_{j=0}^{n}  (-1)^{j} \binom{-\alpha}{j} y_{n-j}.
$
Solving for $y_n$ gives another candidate for a fractional Euler formula:
\begin{equation}
y_n = \frac{1}{(1-h^{\alpha})} \left(y_0 + h^{\alpha} \sum_{j=0}^{n-1} w_{n-j} y_j  \right).
\label{eq:GL:fractional:Euler:limit}
\end{equation}
This is an \textit{implicit} numerical scheme coming from backward differences (revisited later in \eqref{eq:euler:limit:forward:backward}) in the GL construction, so compared to  \eqref{eq:fractional:Euler:limit}, we expect better stability properties from this GL scheme \eqref{eq:GL:fractional:Euler:limit}.
The weights are
$
w_j = (-1)^{j} \binom{-\alpha}{j}.
$
Setting $y_0=1$, \eqref{eq:GL:fractional:Euler:limit} does  satisfy
$
y_n \rightarrow E_{\alpha}(t^{\alpha}),
$
as expected.

\begin{figure*}[ht!]
\centering
\includegraphics[scale=1.1]{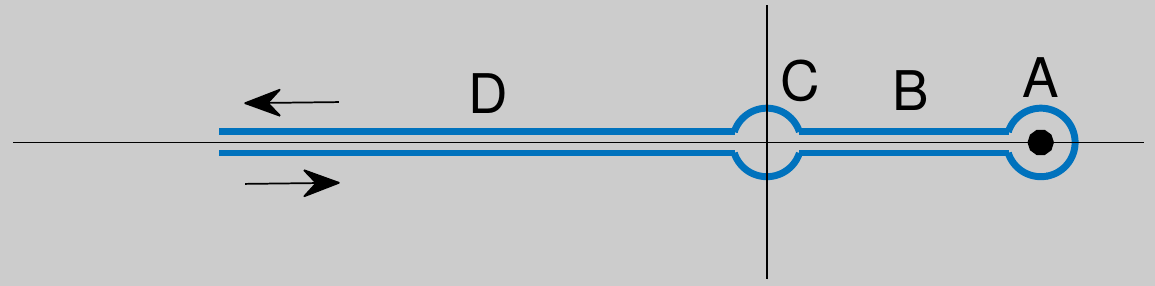}
\caption{Deforming the contour of Mittag-Leffler's representation to collapse the integral to the real axis.
The displayed contour corresponds to the case when the argument to the Mittag-Leffler function is positive, and accompanies the representation of $E_{\alpha}(\bm{+} \lambda t^{\alpha})$ in \eqref{eq:mittag:leffler:positive:argument:integral:little:w}.}
\label{fig:MittagLeffler:contour}
\end{figure*}

\section{\label{sec:Cauchy}A Cauchy integral representation}
 Guiding the search for a discrete construction of a fractional Euler limit has been the  principle that it should limit to a continuous Mittag-Leffler function.
A discretization of the Cauchy integral representation of the Mittag-Leffler function offers another path to this end.

To get from a series representation of the Mittag-Leffler function  \eqref{eq:mittag:leffler} to the Cauchy integral representation, start with a Laplace transform of the series for $E_{\alpha}(\lambda t^{\alpha})$, term by term.
The result is a geometric series that sums to the desired  transform
\begin{equation}
\mathcal{L} \left\{ E_{\alpha}(\lambda t^{\alpha}) \right\} = \frac{s^{\alpha -1}}{ s^{\alpha} - \lambda} =  \frac{1}{s - s^{1-\alpha}\lambda }.
\label{eq:laplace:transform:mittag:leffler:scalar}
\end{equation}
The special case $\lambda=-1$ arises often:
$
\mathcal{L}  \left\{ E_{\alpha}(- t^{\alpha}) \right\} = s^{\alpha -1}/(1+ s^{\alpha}) =  (s^{1-\alpha}+ s)^{-1}.
$
Here the Laplace transform is  $\hat{f}(s) = \mathcal{L} \left\{ f(t) \right\} \equiv \int_0^{\infty} \exp(-st) f(t) \ud t$, and the inverse  transform is $f(t) = \mathcal{L}^{-1} \{ \hat{f}(s) \} \equiv (2 \pi i)^{-1} \int_{\mathcal{C}} \exp(st) \hat{f}(s) \ud s$ where the contour $\mathcal{C}$ is a line parallel to the imaginary axis and to the right of all singularities of $\hat{f}$.
The inverse transform gives
\[
E_{\alpha}(\lambda t^{\alpha}) = \mathcal{L}^{-1} \left\{  \mathcal{L} \left\{ E_{\alpha}(\lambda t^{\alpha}) \right\} \right\} =  \frac{1}{2 \pi i} \int_{\mathcal{C}} e^{st} \frac{s^{\alpha -1}}{ s^{\alpha} - \lambda} \ud s
\]
and after a change of variables this leads to Mittag-Leffler's representation in Table~\ref{tab:exp:mittag:leffler}:
$
 E_{\alpha}(z) = (2 \pi i)^{-1}  \int_{\mathcal{C}} e^{s} s^{\alpha-1} (s^{\alpha}-z)^{-1} \ud s.
$
Here the contour $\mathcal{C}$ must start and end at $- \infty$, and must enclose all singularities and branch points.

  \subsection{\label{sec:Cauchy:branch:cut} Revisiting the Cauchy integral on a branch cut}
Put  $z=\lambda t^{\alpha}$ in Table~\ref{tab:exp:mittag:leffler} to focus on
\begin{equation}
 E_{\alpha}(\lambda t^{\alpha}) = \frac{1}{2 \pi i}  \int_{\mathcal{C}} e^{s} \frac{s^{\alpha-1}}{s^{\alpha}-\lambda t^{\alpha}} \ud s.
\label{eq:Mittag:Leffler:contour:integral:representation:lambda:alpha}
\end{equation}
By collapsing the contour to the real axis (Figure~\ref{fig:MittagLeffler:contour}) we will now arrive at another representation, in \eqref{eq:mittag:leffler:negative:argument:integral} and \eqref{eq:mittag:leffler:positive:argument:integral}.
There are two cases: $E_{\alpha}(\bm{+} \lambda t^{\alpha})$ \textit{is treated separately to} $E_{\alpha}(\bm{-} \lambda t^{\alpha})$.

\subsubsection{A negative argument: $E_{\alpha}( - \lambda t^{\alpha})$}
\begin{figure}[ht!]
\centering
\begin{tabular}{cc}
\includegraphics[scale=0.3]{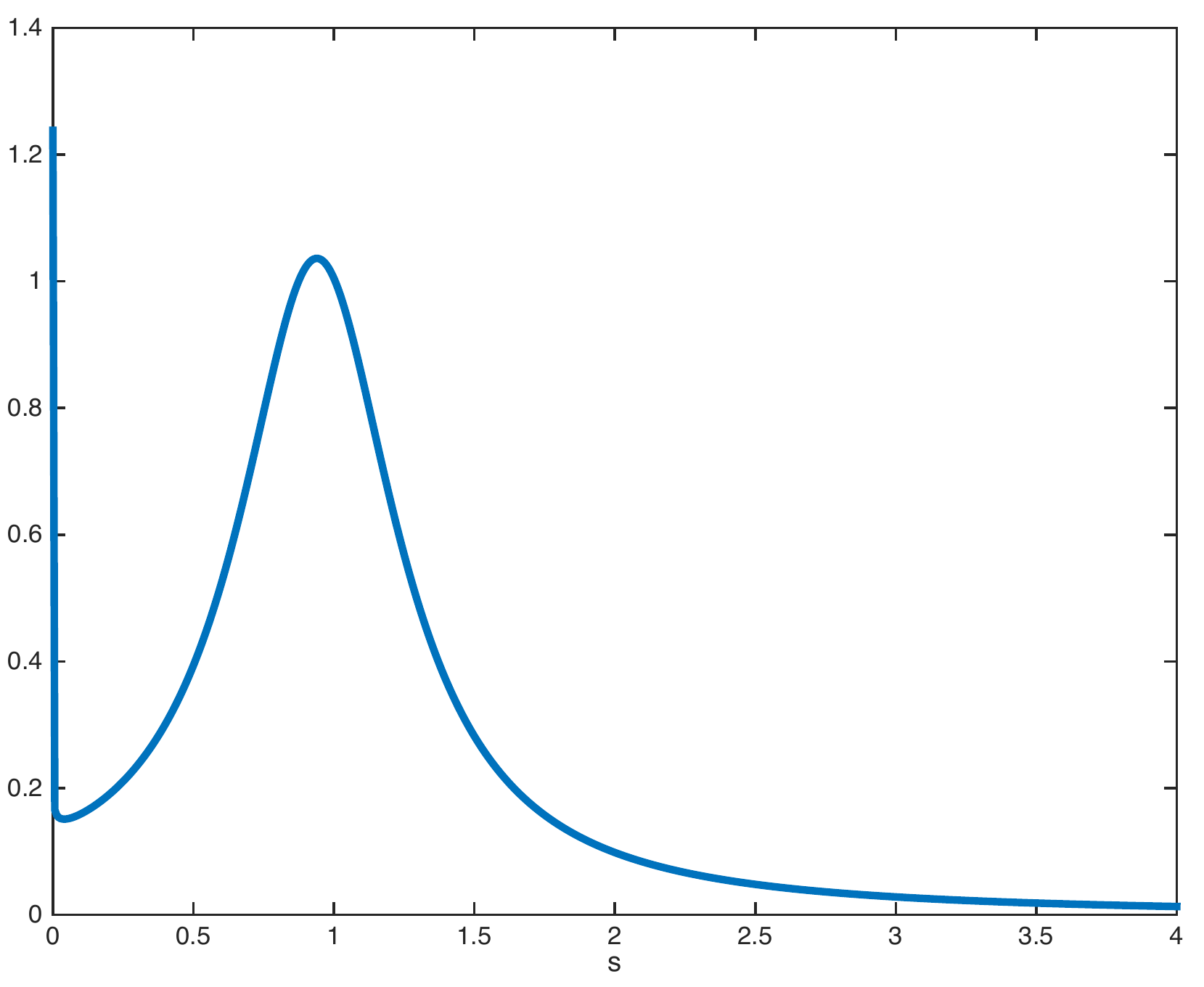} &
\includegraphics[scale=0.3]{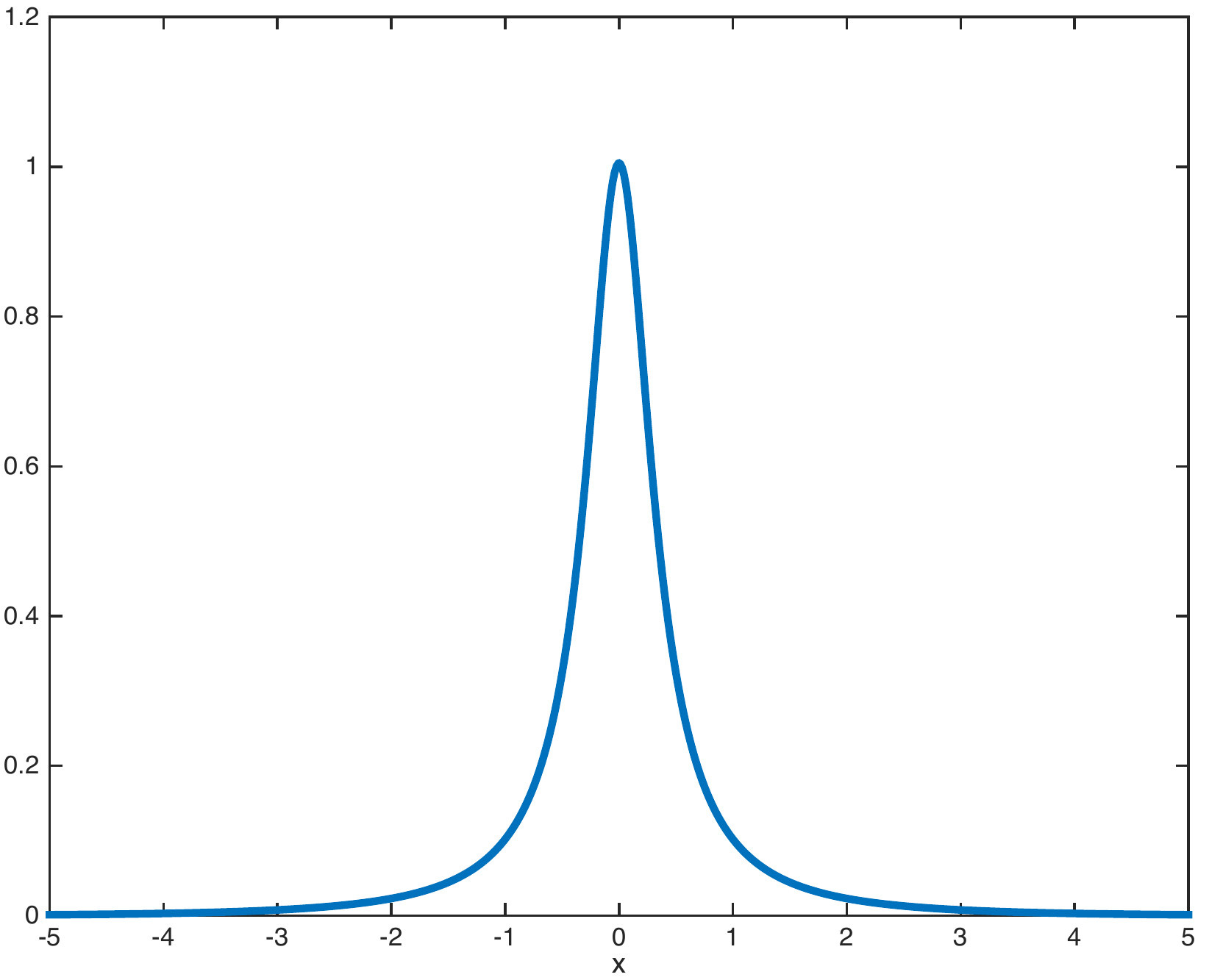}
\end{tabular}
\caption{\textit{Left:} The probability density $w_{\bm{-}}(s)$ of \eqref{eq:mittag:leffler:density:w:minus} appearing in the integral representation of the Mittag-Leffler function with a \textit{negative} argument \eqref{eq:mittag:leffler:negative:argument:integral}.
\textit{Right:} The function $v$ in $w_{\bm{-}}(s) = v(s)/s$, with the substitution  $s = \exp(x)$ as described near equation \eqref{eq:C}.
Parameters: $\alpha=0.9$, $\lambda=1$.
}
\label{fig:MittagLeffler:w:density:minus}
\end{figure}

   This representation is well-known:
 \begin{equation}
E_{\alpha}(- \lambda t^{\alpha}) = \int_0^{\infty}  w_{-}(s)  \exp(-s t) \ud s.
 \label{eq:mittag:leffler:negative:argument:integral}
 \end{equation}
Here $w_{-}(s)$ is the probability density
 \begin{equation}
w_{-}(s) \equiv \lambda \frac{\sin (\alpha \pi)}{\pi} \frac{s^{\alpha-1}}{s^{2 \alpha} + 2 \lambda s^{\alpha} \cos (\alpha \pi) +\lambda^{2}} \ge 0.
 \label{eq:mittag:leffler:density:w:minus}
 \end{equation}
Equation \eqref{eq:mittag:leffler:negative:argument:integral} shows the Mittag-Leffler function as a mixture of exponentials.
An example of the weighting, $w_{-}$, of the components in that mixture is shown in Figure~\ref{fig:MittagLeffler:w:density:minus} (left).
As $\alpha \rightarrow 1$ the weighting converges to Dirac's delta distribution so that the `mixture' becomes a pure exponential with rate $\lambda$, and the Mittag-Leffler function $E_{\alpha}(- \lambda t^{\alpha})$ smoothly transforms to the exponential function $\exp (-\lambda t)$.
From \eqref{eq:mittag:leffler:negative:argument:integral}, $\mathcal{L} \left\{ E_{\alpha}(-\lambda t^{\alpha}) \right\}$ is `the Laplace transform of the Laplace transform of $w_{-}$':
\begin{equation}
  \mathcal{L} \left\{ \mathcal{L} \left\{ w_{-}(s) \right\} \right\} = \mathcal{L} \left\{ E_{\alpha}(-\lambda t^{\alpha}) \right\} = \frac{s^{\alpha-1}}{\lambda +s^{\alpha}}.
\label{eq:laplace:transform:mittag:lefflerminus:lambda}
\end{equation}

\subsubsection{A positive argument: $E_{\alpha}( + \lambda t^{\alpha}) $}

Although \eqref{eq:mittag:leffler:negative:argument:integral} displays two purely real representations, recognizing equality seems to have come historically by passing through the imaginary domain, in line with Hadamard's dictum \cite{LaxZalcmanBook2012}, as we now do again for $E_{\alpha}( + \lambda t^{\alpha}) $.

The contour integral is the sum of the contributions from the four parts marked $A$, $B$, $C$
and $D$ of the deformation of the contour in Figure~\ref{fig:MittagLeffler:contour}.
There is a pole on the real axis at $z = t \lambda^{\frac{1}{\alpha}}$.
Deforming the path near $A$ shows this pole contributes a residue of
$
\exp(t \lambda^{\frac{1}{\alpha}}) / \alpha.
$
Combined, the paths marked by $B$ cancel to make zero net contribution.
The path nearby the origin marked by $C$ can be deformed to approach a circle that makes zero contribution.
The origin is a branch point and the negative real axis is a branch cut.
This is because \eqref{eq:Mittag:Leffler:contour:integral:representation:lambda:alpha} involves a term of the form
\[
s^{\alpha} \equiv \exp \left(\textrm{Log}|s| + i \textrm{Arg}(s) \right),
\]
which is analytic on $\mathbb{C} \backslash (-\infty,0]$.
The two semi-infinite paths, displayed on either side of  the negative real axis, near $D$, do not cancel because of the discontinuity across the negative real axis.
For $s$ on the negative real axis, think of $s=x$ and
$
s^{\alpha} = (-x)^{\alpha} e^{ i \textrm{Arg} (s) \alpha}
$
where $\textrm{Arg} (s) = +\pi$ on the path coming from above the real axis and  $\textrm{Arg} (s) = -\pi$ on the path from below.
The two paths contribute:
\[
\frac{1}{2 \pi i} \int_{\bm{-\infty}}^{\bm{0^-}} \frac{(-x)^{\alpha-1} e^{\bm{-} i \pi (\alpha-1)}}{(-x)^{\alpha} e^{\bm{-} i \pi \alpha} -\lambda t^{\alpha} } e^{x} \ud x \;\;\;
+\;\;\; \frac{1}{2 \pi i} \int_{\bm{0^+}}^{\bm{-\infty}} \frac{(-x)^{\alpha-1} e^{ i \pi (\alpha-1)}}{(-x)^{\alpha} e^{ i \pi \alpha} -\lambda t^{\alpha} } e^{x} \ud x.
\]
The change of variables $x \rightarrow -x$ converts this to an integral on the positive real axis and manipulation brings that integral to the form
$
\bm{-} \int_0^{\infty}  w_{\bm{+}}(s)  \exp(-s t) \ud s
$
with weighting $w_+$ in \eqref{eq:mittag:leffler:density:w:plus}.

Putting all the pieces together gives
 \begin{equation}
E_{\alpha}(\bm{+} \lambda t^{\alpha}) = \frac{\exp(t \lambda^{\frac{1}{\alpha}})}{\alpha} \;\; \bm{-} \int_0^{\infty}  w_{\bm{+}}(s)  \exp(-s t) \ud s
 \label{eq:mittag:leffler:positive:argument:integral:little:w}
 \end{equation}
where the weighting function is:
 \begin{eqnarray}
w_{\bm{+}}(s) &\equiv&  - (- \lambda) \frac{\sin (\alpha \pi)}{\pi} \frac{s^{\alpha-1}}{s^{2 \alpha} + 2 (-\lambda) s^{\alpha} \cos (\alpha \pi) + (-\lambda)^{2}}   \;\; \ge 0.
 \label{eq:mittag:leffler:density:w:plus}
 \end{eqnarray}
The minus signs in \eqref{eq:mittag:leffler:density:w:plus} emphasise that $w_{\bm{+}}(s)$ can be obtained from $w_{\bm{-}}(s)$ in \eqref{eq:mittag:leffler:density:w:minus} by merely changing the sign of $\lambda$, and by including one more overall sign change.
It will be helpful to identify some properties of \eqref{eq:mittag:leffler:positive:argument:integral:little:w}.

Completing the square shows the denominator is positive, and
thus for $0<\alpha<1$, \textit{the weighting function $w_{\bm{+}}(s)$ is positive}.
For moderately large $t$, the dominant contribution to $E_{\alpha}(+ \lambda t^{\alpha}) $ in  \eqref{eq:mittag:leffler:positive:argument:integral:little:w} comes from $\exp(t \lambda^{\frac{1}{\alpha}}) / \alpha$, while
 the branch cut integral is the smaller moiety.
As $\alpha \rightarrow 1$, the Mittag-Leffler function tends to the exponential so the contribution from the integral in
\eqref{eq:mittag:leffler:positive:argument:integral:little:w}  tends to zero.

When numerically integrating either density $w_{\bm{-}}(s)$ or $w_{\bm{+}}(s)$, a challenge may arise due to the the scaling as  $s \rightarrow 0^+$, and due to slow decay as $s \rightarrow \infty$ (Figure~\ref{fig:MittagLeffler:w:density:minus}, left).
 A possible remedy is to write $w_{\bm{+}}(s) = v(s)/s$, and then to use the substitution $s = \exp(x)$.
  In the new variables the scalings go exponentially to zero   (Figure~\ref{fig:MittagLeffler:w:density:minus}, right): $v(\exp(x)) \sim \exp(\alpha x)$ as  $x \rightarrow - \infty$ and $v(\exp(x)) \sim \exp(-\alpha x)$ as $x \rightarrow \infty$.
  For example, we can numerically evaluate  the integral on the right  of
   \begin{equation}
  0 \le C \equiv \int_0^{\infty} w_+(s) ds = \int_0^{\infty} \frac{v(s)}{s} \, ds = \int_{-\infty}^{\infty} v(\exp(x)) dx ,
   \label{eq:C}
   \end{equation}
which can reasonably be handled by first truncating to a finite subset of the real line and then using an equally-spaced quadrature rule such as the trapezoidal rule.

It transpires that
\begin{equation}
C = \frac{1}{\alpha} - 1
\label{eq:normalization:constant:C}
\end{equation}
 is a normalisation constant for the probability density
 \begin{equation}
W_{\bm{+}}(s)  \equiv   \frac{1}{C} w_{\bm{+}}(s)  =  \frac{\alpha}{1-\alpha}w_{\bm{+}}(s).
\label{eq:mittag:leffler:density:w:plus:big:W}
 \end{equation}
For example, if $\alpha=1/2$ then $C=1$, and comparing \eqref{eq:mittag:leffler:density:w:minus} and \eqref{eq:mittag:leffler:density:w:plus}, this is a special case where the densities are the same: $W_+ = w_{\bm{+}} = w_{-}$.

\begin{figure}[ht!]
\centering
\includegraphics[scale=0.5]{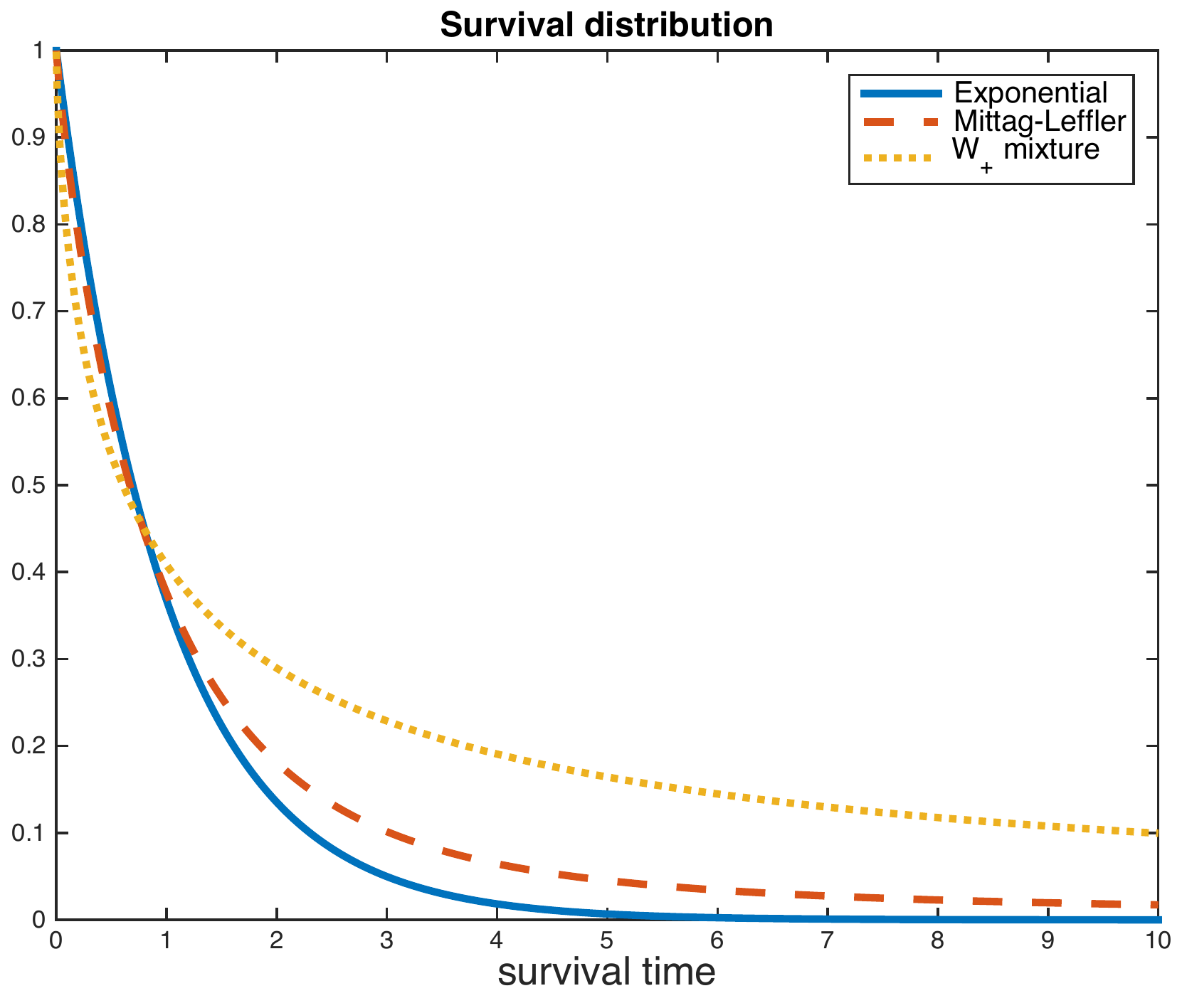}
\caption{Comparison of survival times from three distributions.
$(i)$ A pure, unit rate exponential distribution.
$(ii)$ A Mittag-Leffler distribution, which is the mixture of exponentials in $\int_0^{\infty}  w_{-}(s)  \exp(-s t) \ud s$ of  \eqref{eq:mittag:leffler:negative:argument:integral}
with the probability density $w_{-}(s)$  of  \eqref{eq:mittag:leffler:density:w:minus}.
$(iii)$ The mixture of exponentials $ \phi_{W_{+}}(t) = \int_0^{\infty}  W_{\bm{+}}(s)  \exp(-s t) \ud s$ in \eqref{eq:W:plus:exponential:mixture}  with density $W_{+}$  as in  \eqref{eq:mittag:leffler:density:w:plus:big:W}.
Parameters: $\alpha=0.9$, $\lambda=1$.
}
\label{fig:waiting:time:compare}
\end{figure}

The representation of the Mittag-Leffler function \eqref{eq:mittag:leffler:positive:argument:integral:little:w} becomes
 \begin{equation}
E_{\alpha}(\bm{+} \lambda t^{\alpha}) = \frac{\exp(t \lambda^{\frac{1}{\alpha}})}{\alpha} \;\; \bm{-} \frac{1-\alpha}{\alpha}  \int_0^{\infty}  W_{\bm{+}}(s)  \exp(-s t) \ud s.
 \label{eq:mittag:leffler:positive:argument:integral}
 \end{equation}
 The last term in \eqref{eq:mittag:leffler:positive:argument:integral}  is a mixture of exponential distributions
 \begin{equation}
 \phi_{W_{+}}(t) \equiv  \int_0^{\infty}  W_{\bm{+}}(s)  \exp(-s t)\ud s.
  \label{eq:W:plus:exponential:mixture}
 \end{equation}
 Compared to a pure exponential distribution, this $W_{+}$-mixture has heavier tails, as does the Mittag-Leffler distribution (Figure~\ref{fig:waiting:time:compare}).
Thus, collapsing the Cauchy integral to the real axis leads naturally to a density $W_{+}$ associated with a Mittag-Leffler function of a positive argument,  analogous to the way that $w_{-}$  arises from the Cauchy integral associated with a negative argument.
Connections to \textit{Wright functions} and to \textit{Fox H-functions} are important but not discussed here \cite{MainardiBook2010,FoxHFunctionBook2009}.
Next, we derive formulas to facilitate sampling from the $W_+$ density.

\paragraph{The cumulative distribution function}
 Denote the desired integral
\begin{equation}
g(s) \equiv \int w_{\bm{+}}(s) \ud s =
    \int - (- \lambda) \frac{\sin (\alpha \pi)}{\pi} \frac{s^{\alpha-1}}{s^{2 \alpha} + 2 (-\lambda) s^{\alpha} \cos (\alpha \pi) + (-\lambda)^{2}}  \ud s.
\label{eq:w:plus:integral}
 \end{equation}
There are two cases for  $g$, depending on whether or not $\alpha$ is bigger than $1/3$.
If $0<\alpha<1/3$ then
\begin{equation}
g(s) = \int w_{\bm{+}}(s) ds =
\begin{cases}
\tilde{g}(s) -  \frac{1}{ \alpha} & \mbox{if }  0 < s < s_1 \\
\tilde{g}(s) - \frac{1}{2 \alpha} & \mbox{if }  s_1 < s < s_2 \\
\tilde{g}(s)  & \mbox{otherwise }  (s_2< s).
\end{cases}
\label{eq:ArcTan:w:plus:exact:integral:correction:case:1}
\end{equation}
Otherwise, if $1/3<\alpha<1$ then
\begin{equation}
g(s) = \int w_{\bm{+}}(s) ds =
\begin{cases}
\tilde{g}(s) - \frac{1}{2 \alpha} & \mbox{if }  0 < s < s_2 \\
\tilde{g}(s)  & \mbox{otherwise }  (s_2< s).
\end{cases}
\label{eq:ArcTan:w:plus:exact:integral:correction:case:2}
\end{equation}
Here $s_1$ and $s_2$ are defined in \eqref{eq:s:1} and \eqref{eq:s:2},
\begin{equation}
  s_1 \equiv (2 \lambda \cos (\alpha \pi ) - \lambda)^{\frac{1}{\alpha}},
 \label{eq:s:1}
\end{equation}
\begin{equation}
s_2 \equiv \lambda^{\frac{1}{\alpha}},
 \label{eq:s:2}
\end{equation}
and $\tilde{g}$ is as in \eqref{eq:ArcTan:w:plus:exact:integral:mathematica:f}\footnote{\texttt{Mathematica} evaluates \eqref{eq:w:plus:integral}  as \eqref{eq:ArcTan:w:plus:exact:integral:mathematica:f}, but for our purpose
  the correct integral is \eqref{eq:ArcTan:w:plus:exact:integral:correction:case:1} or  \eqref{eq:ArcTan:w:plus:exact:integral:correction:case:2}. }:
 \begin{equation}
 \tilde{g}(s) \equiv
 \frac{\textrm{ArcTan}\left( \frac{-\lambda + s^\alpha - 2 \lambda \cos(\alpha \pi)  }{
  \lambda + s^\alpha - 2 \lambda \cos (\alpha \pi )} \tan (\alpha \pi/2 ) \right)  }{  2 \alpha \pi}
\; -\frac{        \textrm{ArcTan}\left(\frac{\lambda + s^\alpha}{-\lambda + s^\alpha}  \tan(\alpha \pi/2) \right) }{  2 \alpha \pi}.
\label{eq:ArcTan:w:plus:exact:integral:mathematica:f}
\end{equation}
Notice $s_1<s_2.$
Usually, a function defined by separate cases is not continuous, so it might be surprising to learn that  \textit{our function $g$, defined in \eqref{eq:ArcTan:w:plus:exact:integral:correction:case:1} or \eqref{eq:ArcTan:w:plus:exact:integral:correction:case:2}, is continuous!}
We are allowing $\textrm{ArcTan}(\pm \infty)= \pm \pi/2$ at the singularities $s_1$ and $s_2$.
Note also that $g(s)<0$, and that $g$ is increasing, with $g \rightarrow 0$ as $s \rightarrow \infty$.
The normalisation constant $C=\int_0^{\infty} w_+(s) ds = g(\infty) - g(0) = 0 - (1-\frac{1}{\alpha})$ can be found by evaluating  \eqref{eq:ArcTan:w:plus:exact:integral:correction:case:1} or  \eqref{eq:ArcTan:w:plus:exact:integral:correction:case:2} at the two limits, as $s \rightarrow 0^+$ and $s \rightarrow \infty$.
Observe that $C \rightarrow 0$  as $\alpha \rightarrow 1$, and
at the other extreme, $C \rightarrow \infty$ as $\alpha \rightarrow 0$.
The larger $C$, the `more fractional' the calculus.

At last, with $C$ in \eqref{eq:normalization:constant:C} define  $G:(0, \infty) \rightarrow (0,1)$ by
\begin{equation}
G(T) \equiv \int_0^T W_{+}(s) \ud s = \frac{g(T) - g(0)}{C} = 1 + \frac{g(T)}{C}.
\label{eq:cumulative:distribution:function:G}
\end{equation}
so that $G$ is the cumulative distribution function associated with $W_{+}$ of \eqref{eq:mittag:leffler:density:w:plus:big:W}.
(Here $g$ is given by \eqref{eq:ArcTan:w:plus:exact:integral:correction:case:1} if $0<\alpha<1/3$, or by  \eqref{eq:ArcTan:w:plus:exact:integral:correction:case:2} if $1/3<\alpha<1$.)

\paragraph{The \textit{inverse} cumulative distribution function}
As usual, the inverse function $G^{-1}: (0,1) \rightarrow (0, \infty)$, is defined by the property that, when $u=G(T)$, we have $T=G^{-1}(u)$.
The inverse is the positive real root
 \begin{equation}
 T = G^{-1}(u) = \left( Q(C u) \right)^{\frac{1}{\alpha}} = Q^{\frac{1}{\alpha}}
\label{eq:cumulative:distribution:inverse:function:G}
\end{equation}
where $Q: (0,C) \rightarrow (0, \infty)$ is given by
\begin{eqnarray}
& &Q(v) = \lambda \frac{     \sin (\alpha \pi (1 + 2v)) -  \sin(\alpha \pi) }{\sin \left( 2 \alpha \pi (1 + v) \right)}.
\label{eq:cumulative:distribution:inverse:function:Q}
\end{eqnarray}

\paragraph{A Laplace transform associated with $W_{+}$}
Let \eqref{eq:W:plus:exponential:mixture} define a survival distribution.
The Laplace transform of this survival function is the Laplace transform of the Laplace transform of $W_+ $: \;
$
\hat{\phi}_{W_{+}}(t) = \mathcal{L} \left\{ \phi_{W_{+}}(t) \right\} =  \mathcal{L} \left\{ \mathcal{L} \left\{ W_{+}(s) \right\} \right\}.
$
By \eqref{eq:mittag:leffler:positive:argument:integral},
$
\phi_{W_{+}}(t)  =  \alpha (1 - \alpha)^{-1} ( \exp(t \lambda^{\frac{1}{\alpha}}) / \alpha -  E_{\alpha}(\bm{+} \lambda t^{\alpha}) )
$
so $\mathcal{L} \left\{ \phi_{W_{+}}(t) \right\} = $
\begin{equation}
   \frac{\alpha}{1 - \alpha} \left( \frac{1}{\alpha}  \mathcal{L} \left\{ \exp(t \lambda^{\frac{1}{\alpha}}) \right\}  -  \mathcal{L} \left\{E_{\alpha}(\bm{+} \lambda t^{\alpha}) \right\} \right)
= \frac{\alpha}{1 - \alpha}  \left( \frac{1}{\alpha (-\lambda^{\frac{1}{\alpha}}+s)}  - \frac{s^{\alpha-1}}{-\lambda + s^{\alpha}} \right).
\label{eq:laplace:transform:W:plus}
\end{equation}
Unlike for the exponential function or for the Mittag-Leffler function, the authors do not know if the corresponding ratio of transforms (significant later in \eqref{eq:memory:function}), $\hat{K}_{W_{+}} (s) = \hat{\phi}_{W_{+}}(s)/\hat{\psi}_{W_{+}}(s)  =  \hat{\phi}_{W_{+}}(s)/(s \hat{\phi}_{W_{+}}(s) -1)$, has a simple interpretation.

\subsection{A fractional Euler limit, again}
Inspired by the representation of the continuous Mittag-Leffler function as a mixture of exponentials in \eqref{eq:mittag:leffler:negative:argument:integral} for a negative argument, or \eqref{eq:mittag:leffler:positive:argument:integral} for a positive argument, it is  natural to ask the question,  is it possible to write the discrete fractional limit in the form of a weighted sum of regular Euler limits?
The answer is `yes' and here are two examples.

\textit{Fractional decay, again.}
Instead of \eqref{eq:fractional:Euler:limit:minus}, we now propose a different discrete fractional generalisation of the Euler formula, namely
\begin{equation}
\sum_k w_k (1- s_k t/n)^n .
\label{eq:weighted:sum:euler:limits:negative}
\end{equation}

\textit{Fractional growth, again.}
Likewise, instead of \eqref{eq:fractional:Euler:limit}, we also now propose a fractional generalization of the Euler formula in the case of a positive argument, namely
\begin{equation}
  \frac{\left(1+ (t \lambda^{\frac{1}{\alpha}})/n \right)^n}{\alpha} \; \bm{-} \frac{1-\alpha}{\alpha}   \sum_k W_k (1-s_k t/n)^n .
\label{eq:weighted:sum:euler:limits:positive}
\end{equation}

In both  cases of growth and of decay, it is desirable that the proposed sum converges to the corresponding integral  for large $n$.
That integral  is \eqref{eq:mittag:leffler:negative:argument:integral} in the case of decay, or \eqref{eq:mittag:leffler:positive:argument:integral} in the case of growth.
That is, it is desirable to have both the following discrete-to-continuous limits,  as $n \rightarrow \infty$:
$(i)$ \,  $w_k \longrightarrow w(s)$ \, and \;
$(ii)$\, $ \sum_k w_k (1 \bm{-} s_k t/n)^n  \longrightarrow  \int_0^{\infty}  w(s)  \exp(- s t) \ud s$. \,
Here $k$ would depend on $n$, and $w=w_{-}$ of \eqref{eq:mittag:leffler:density:w:minus} in the case of decay, or $w=W_{+}$ of \eqref{eq:mittag:leffler:density:w:plus:big:W} in the case of growth.
Always, the discrete weights $w_k = w(s_k)$ integrate to one, i.e. $w_k>0$ and $1 = \sum_k w_k (s_k - s_{k-1})$.
We have merely suggested the general form in \eqref{eq:weighted:sum:euler:limits:negative} and \eqref{eq:weighted:sum:euler:limits:positive} --- these formulations can be interpreted as quadrature rules applied to corresponding integrals so there remains the potential for myriad variations in the details not specified here, such as the spacing and the number of grid points.
These would not usually be good numerical schemes for the reasons highlighted in Figure~\ref{fig:MittagLeffler:w:density:minus}.
However, as a weighted sum of regular Euler limits, this form is an especially satisfying fractional generalization.

\subsection{\label{sec:completely:monotone} Complete monotonicity}
A smooth  function $f$ is \textit{completely monotone}  if all derivatives are monotone: $(-1)^n f^{(n)} \ge 0$. 
Exponential decay, $e^{-t}$, is the prototypical example.
Bernstein's theorem tells us that all completely monotone functions are representable as a mixture\footnote{A mixture usually refers to a finite sum, or at most countably infinite sum, whereas here in our integral with an exponential kernel we are  allowing an uncountably infinite `mixture'.} of exponentials; being completely monotone is equivalent to being the real Laplace transform of a non-negative function.
For the Mittag-Leffler function with a negative argument, this non-negative function $w$ is explicitly known in \eqref{eq:mittag:leffler:negative:argument:integral} and \eqref{eq:mittag:leffler:density:w:minus}.

The Mittag-Leffler function with a positive argument does not have the completely monotone property.
However, the weighting functions $w_{+}(s)$ in \eqref{eq:mittag:leffler:density:w:plus} or $W_+$ in \eqref{eq:mittag:leffler:density:w:plus:big:W} are positive so
$
 \phi_{W_+} = \int_0^{\infty}  W_{\bm{+}}(s)  \exp(-s t) \ud s,
$
 the transform of $W_+$  in  \eqref{eq:W:plus:exponential:mixture}, does have the property.
Thus \eqref{eq:mittag:leffler:positive:argument:integral} shows us the Mittag-Leffler function of a positive argument as a combination of two functions, one of which is completely monotone.

Both the exponential and the Mittag-Leffler function are entire functions, and both are completely monotone in the case of a negative argument, and both lose this monotone property in the case of a positive argument.
Amongst our candidates for the fractional generalisation of the Euler limit,  \eqref{eq:weighted:sum:euler:limits:negative} is more amenable to  emphasising such shared properties of the exponential and the Mittag-Leffler function.
As an example, suppose $t>0$ and $\lambda>0$ so for all sufficiently large $n$,  $0<(1 \bm{-} \lambda t/n)<1$.
A very useful application of the Euler limit formula
$
(1\bm{-} \lambda t/n)^n \longrightarrow \exp( \bm{-} \lambda t),
$
is to make clear that the exponential of a negative argument satisfies
 $
0<\exp(\bm{-} \lambda t) <1.
 $
The Mittag-Leffler function shares this property:
\[
0< E_{\alpha}(\bm{-} \lambda t^{\alpha}) <1.
\]
Representations of the Mittag-Leffler function as a weighted integral of exponentials where the weight is a nonnegative density  as in  \eqref{eq:mittag:leffler:negative:argument:integral}, or the discrete fractional Euler limit in \eqref{eq:weighted:sum:euler:limits:negative},  make it clear that the Mittag-Leffler function does indeed have this property.
In contrast, although the fractional generalisations in \eqref{eq:fractional:Euler:limit} or \eqref{eq:GL:fractional:Euler:limit} have their own virtues, it is nearly impossible to discern this important property from those formulations.
This makes a discretization such as suggested by \eqref{eq:weighted:sum:euler:limits:negative}, a more attractive fractional generalisation.

 \section{\label{sec:probabilistic} A probabilistic interpretation}

\subsection{Sampling from Mittag-Leffler distributions or the $W_+$ mixture}
A Mittag-Leffler distribution is characterised by three functions.
The \textit{survival time}  is
\begin{equation}
\phi(t) = E_{\alpha}(- \lambda t^{\alpha})
\label{eq:survival:distribution:mittag:leffler}
\end{equation}
so  the  \textit{waiting time} density ($\psi(t)=  - \ud \phi(t) / \ud t = - \ud E_{\alpha}(- \lambda t^{\alpha})  / \ud t $) is
\begin{equation}
\psi(t) = \lambda t^{\alpha -1} E_{\alpha,\alpha}(- \lambda t^{\alpha})
\label{eq:waiting:time:mittag:leffler}
\end{equation}
where $E_{\alpha,\beta}(z) \equiv \sum_{k=0}^{\infty} \frac{z^{k}}{\Gamma(\alpha k +\beta)} $ is the two-parameter Mittag-Leffler function, and the \textit{cumulative distribution} function ($cdf(t) = 1- \phi(t)  = \int_0^t \psi(s) \ud s$) is
\begin{equation}
cdf(t) = 1 - E_{\alpha}(- \lambda t^{\alpha}) .
\label{eq:cumulative:distribution:mittag:leffler}
\end{equation}
All three functions are nonnegative and, as a probability density  $1 = \int_0^{\infty} \psi(t) \ud t $.
Putting $\alpha \leftarrow 1$ in these formulas recovers the survival time distribution, $cdf$,  and waiting time density, corresponding to an exponential distribution with parameter $\lambda$ and mean value $1/\lambda$.
Unlike the exponential case, equation \eqref{eq:waiting:time:mittag:leffler} shows that the Mittag-Leffler function and hence also the solution of \eqref{eq:fractional:ode} is not differentiable at $t=0^+$.
In general care must be taken when differentiating near zero, as happens later in \eqref{eq:p:derivative}.

The survival time  \eqref{eq:survival:distribution:mittag:leffler} is also a mixture of exponentials, $\phi(t) = E_{\alpha}(- \lambda t^{\alpha})= \int_0^{\infty}  w_{-}(s)  \exp(-s t) \ud s$  (\eqref{eq:mittag:leffler:negative:argument:integral},  \eqref{eq:mittag:leffler:density:w:minus}).
Such mixtures are amenable to methods for fast simulation of geometric stable distributions involving products of independent random variables  \cite{BelinskiyKozubowski2000,KozubowskiRachev1999,Fulger2008,Kozubowski2000,LaplaceDistributionBook}.
Let $F(T) \equiv \int_0^T w_{-}(s) ds$ denote the cumulative distribution function of $w_{-}(s)$.
Given $u \in (0,1)$ the inverse function, that solves $F(T)=u$ for $T$, is known to be
$
T = F^{-1}(u) = \lambda^\frac{1}{\alpha}  ( \sin(\pi \alpha)/\tan(\pi \alpha (1-u)) - \cos(\pi \alpha) )^\frac{1}{\alpha}.
$
This permits the known fast sampling procedure for the Mittag-Leffler distribution, via
\begin{equation}
\tau \; \sim \;\; - \left(\frac{1}{\lambda}\right)^\frac{1}{\alpha}  \left( \frac{\sin(\pi \alpha)}{\tan(\pi \alpha (1-u_1))} - \cos(\pi \alpha) \right)^\frac{1}{\alpha}  \log (u_2).
\label{eq:sampling:mittag:leffler}
\end{equation}
Here $u_1$ and $u_2$ are independent samples from the uniform distribution on $(0,1)$.
As $\alpha \rightarrow 1$, \eqref{eq:sampling:mittag:leffler} reduces to the familiar formula $\tau \sim - \log (u)/ \lambda$ for sampling from an exponential distribution with density $\lambda \exp(-\lambda t)$.

One way to understand \eqref{eq:sampling:mittag:leffler} is as the product of two independent random variables, $Z = XY$, with density  $f_Z(z) = \int_{-\infty}^{\infty} f_X(x)  f_Y(\frac{z}{x}) |x|^{-1} \ud x $.
Now read \eqref{eq:sampling:mittag:leffler} as $XY$ where $X \sim ( \sin(\pi \alpha)/\tan(\pi \alpha (1-u_1)) - \cos(\pi \alpha) )^\frac{1}{\alpha} $ is the inverse transform method for sampling from a density $f_X(x) = w_{-,1}(x)$, and where $w_{-,1}$ is  \eqref{eq:mittag:leffler:density:w:minus:1} (the special case of the density  with $\lambda=1$ in \eqref{eq:mittag:leffler:density:w:minus}, but note that this does not assume $\lambda=1$ in  \eqref{eq:sampling:mittag:leffler}) and where $ Y \sim  -  \log (u_2) / \lambda^{\frac{1}{\alpha}}$
 is the familiar inverse transform method for sampling from an exponential distribution with density $f_Y = \lambda^{\frac{1}{\alpha}} \exp(-\lambda^{\frac{1}{\alpha}} y)$.
 Then $f_Z(z) = \int_{0}^{\infty} w_{-,1}(x) \lambda^{\frac{1}{\alpha}} \exp(-\lambda^{\frac{1}{\alpha}}\frac{z}{x}) |x|^{-1} \ud x $.
With the change of variables  $s = 1/x$ this becomes $f_Z(z) = \int_{\infty}^{0} w_{-,1}(1/s) \lambda^{\frac{1}{\alpha}}  \exp(- \lambda^{\frac{1}{\alpha}} s z) s  (-s^2) \ud s $.
Noticing   $w_{-}(s) = s^2 w_{-}(1/s)$, and replacing $z$ by $t$ this becomes
\begin{equation}
 \int_{0}^{\infty} w_{-,1}(s)  s \lambda^{\frac{1}{\alpha}}  \exp(-s\lambda^{\frac{1}{\alpha}}t)  \ud s .
\label{eq:ML:exponential:mixture:sum:form}
\end{equation}
This is   $- \ud/ \ud t  \,   \int_{0}^{\infty} w_{-,1}(s)    \exp(-s\lambda^{\frac{1}{\alpha}}t)  \ud s  =  - \ud E_{\alpha}(- \lambda t^{\alpha})  / \ud t =- \ud \phi(t) / \ud t$, which is the derivative of the representation of $\phi(t)$ in \eqref{eq:w:1:exp:lambda:1:alpha:representation}.
Thus \eqref{eq:ML:exponential:mixture:sum:form} is another representation of $\psi(t)$ because the waiting time density of  \eqref{eq:waiting:time:mittag:leffler} is always $- \ud \phi(t) / \ud t $.
This confirms  \eqref{eq:sampling:mittag:leffler} does indeed sample Mittag-Leffler waiting times.


Instead of the above product form, we could think of the Mittag-Leffler density \eqref{eq:ML:exponential:mixture:sum:form} as a sum of densities of exponential random variables.
This suggests a recipe: sample the exponential random variable with parameter $s\lambda^{\frac{1}{\alpha}}$, where $s$ is sampled according to the density $w_{-,1}(s)$.
The recipe could instead just as well sample the exponential random variable with the parameter in which $s$ is replaced by $1/s$.
This is because the same change of variables shows $\psi(t) = \int_{0}^{\infty} w_{-,1}(s)  s \lambda^{\frac{1}{\alpha}}  \exp(-s\lambda^{\frac{1}{\alpha}}t)  \ud s = \int_{0}^{\infty} w_{-,1}(s)   \lambda^{\frac{1}{\alpha}} s^{-1}  \exp(-\lambda^{\frac{1}{\alpha}} s^{-1} t)  \ud s.$
The sampling formula displayed in \eqref{eq:sampling:mittag:leffler} corresponds to the latter choice.

This same pattern that works for Mittag-Leffler also works for the survival function
$
\phi_{W_{+}}(t) =  \int_0^{\infty}  W_{\bm{+}}(s)  \exp(-s t)\ud s
$
(Figure~\ref{fig:waiting:time:compare}, \eqref{eq:W:plus:exponential:mixture}), which is also a mixture of exponentials.
Fortunately, we identified both the cumulative distribution of $W_{+}$  in \eqref{eq:cumulative:distribution:function:G} and its inverse in \eqref{eq:cumulative:distribution:inverse:function:G} and \eqref{eq:cumulative:distribution:inverse:function:Q}.
Thanks to \eqref{eq:cumulative:distribution:inverse:function:Q}  we again have a fast Monte Carlo procedure to sample waiting times from the $W_+$ mixture:
$ C \leftarrow \frac{1}{\alpha} -1$, $v \leftarrow  C u_1$, and
\begin{equation}
\tau \sim  \;\; \; - \left(\frac{1}{\lambda}\right)^\frac{1}{\alpha}  \left(   \frac{     \sin (\alpha \pi (1 + 2v)) -  \sin(\alpha \pi) }{\sin \left( 2 \alpha \pi (1 + v) \right)}   \right)^\frac{1}{\alpha}  \log (u_2)
\label{eq:sampling:W:plus:mixture:of:exponentials}
\end{equation}
where $u_1,u_2$ are independent uniform random variables on $(0,1)$.

\subsection{Matrix arguments and the special density when $\lambda=1$}
The notation $w_{-}$ suppresses the dependence of the density on the parameters $\alpha$ and $\lambda$.
Put $\lambda=1$ in \eqref{eq:mittag:leffler:negative:argument:integral} to obtain
\begin{equation}
E_{\alpha}(- t^{\alpha}) = \int_0^{\infty}  w_{-,1}(s)  \exp(-s t) \ud s
\label{eq:mittag:leffler:negative:argument:integral:1}
 \end{equation}
where  $w_{-,1}(s)$ is the special case of the probability density $w_{-}(s)$ when $\lambda=1$ in \eqref{eq:mittag:leffler:density:w:minus}:
 \begin{equation}
w_{-,1}(s) \equiv  \frac{\sin (\alpha \pi)}{\pi} \frac{s^{\alpha-1}}{s^{2 \alpha} + 2  s^{\alpha} \cos (\alpha \pi) + 1} \ge 0.
 \label{eq:mittag:leffler:density:w:minus:1}
 \end{equation}
Make the replacement $t \leftarrow \lambda^{\frac{1}{\alpha}} t$ (now allowing $\lambda$ to take any positive value) to obtain
 \begin{equation}
 E_{\alpha}( - \lambda t^{\alpha}) = \int_0^{\infty}  w_{-,1}(s)  \exp(-s \lambda^{\frac{1}{\alpha}}  t) \ud s.
 \label{eq:w:1:exp:lambda:1:alpha:representation}
 \end{equation}
 Alternatively, we could arrive at  \eqref{eq:w:1:exp:lambda:1:alpha:representation} by the change of variables $s \leftarrow \lambda^{\frac{1}{\alpha}}s$ in  \eqref{eq:mittag:leffler:negative:argument:integral}, and by noticing  $ \lambda^{\frac{1}{\alpha}} w_{-} (\lambda^{\frac{1}{\alpha}} s) = w_{-,1}(s)$.
Compared to the representation of  $E_{\alpha}( - \lambda t^{\alpha})$  in \eqref{eq:mittag:leffler:negative:argument:integral}, here in \eqref{eq:w:1:exp:lambda:1:alpha:representation} the weighting density $w_{-,1}(s)$ does \textit{not} depend on $\lambda$ (but the exponential $\exp(st)$ in \eqref{eq:mittag:leffler:negative:argument:integral} has been replaced  by something that does, $\exp(-s \lambda^{\frac{1}{\alpha}}  t)$, so the dependency has been `moved' from the weighting density to the exponential).

\textit{A Mittag-Leffler matrix function.}
There are various issues when defining a matrix function \cite{NicholasHighamBook08}.
 When defining a Mittag-Leffler function of a matrix $\mathbb{A} \in \mathbb{R}^{n \times n}$, one possibility is a series that replaces $z \in \mathbb{C}$ with $\mathbb{A}$ in \eqref{eq:mittag:leffler}, while another possibility is to apply the scalar Mittag-Leffler function to eigenvalues together with a diagonalization.
Alternatively, if the eigenvalues of the matrix are positive (recall that  the exponential mixture representations \eqref{eq:mittag:leffler:negative:argument:integral} and \eqref{eq:mittag:leffler:positive:argument:integral} coming from the Cauchy integral  have different forms depending on the sign of the argument) then it will usually be meaningful to make the replacement $\lambda \leftarrow \mathbb{A}$ in \eqref{eq:w:1:exp:lambda:1:alpha:representation}  to express the Mittag-Leffler function of a matrix as a weighted sum of exponentials of a power of that matrix.

\textit{A graph Laplacian matrix.}
Always in this article a matrix $\mathbb{A}$ has two properties:
\begin{enumerate}[i]
\item  \; Off-diagonal entries that are positive or zero: \;\;\;\;\;\;\;\;\; $a_{ij} \ge 0 \;\;\; (i \ne j)$
\item  \; Diagonal entries  that ensure columns sum to zero: \; $a_{jj} = - \sum_{i, i \ne j} a_{ij}$.
\end{enumerate}
When symmetric, this matrix is the familiar \textit{graph Laplacian}.
 Many authors use this same terminology for both symmetric and nonsymmetric cases.
The spectra of such matrices is entirely in the left-half complex plane, and they are examples of sectorial operators.
We often further assume the matrix $\mathbb{A}$ has distinct, negative and real eigenvalues, and one zero eigenvalue.
With this notation\footnote{A numerical analyst will be frustrated, preferring the opposite sign convention.
The reason for this choice of signs is that it is common in the Markov process literature to denote an infinitesimal generator this way, although the usual `$Q$-matrix' is the transpose of $\mathbb{A}$.} it is the matrix $(-\mathbb{A})$ that we think of as `positive' (!)  so the representation is
 \begin{equation}
E_{\alpha}(\mathbb{A}t^{\alpha}) = \int_0^{\infty}  w_{-,1}(s)  \exp(-s (-\mathbb{A})^{\frac{1}{\alpha}}  t) \ud s.
\label{eq:mittag:leffler:negative:argument:integral:1:matrix:argument}
 \end{equation}
This representation \eqref{eq:mittag:leffler:negative:argument:integral:1:matrix:argument} emphasises the role of $(-\mathbb{A})^{\frac{1}{\alpha}}$, and there are important connections to \textit{subordination} (allowing a solution for a particular $\alpha$ to be represented as an integral of the solution corresponding to a different $\alpha$, for example, and related to the monotone properties of section \ref{sec:completely:monotone}) and to  powers of operators not discussed here \cite{Kat76,YosidaFunctionalAnalysisBook,AbadiasMianaSubordination2015,BaeumerMeerschaertSubordination,GorenfloMainardiVivoli2009,BochnerBook}, \cite[Ch. 4]{PrussBook2012}.

When $\alpha=1$,  $w_{-,1}$ is the Dirac distribution centred at $1$ so that \eqref{eq:mittag:leffler:negative:argument:integral:1:matrix:argument} is the familiar matrix exponential $\exp(\mathbb{A}t)$.
 For this class of matrices, $\exp(\mathbb{A}t)$ is a \textit{stochastic matrix} (columns are probability vectors) associated with a Markov process on discrete states in continuous time.
To see that the entries of the matrix exponential are all nonnegative, we could first examine the Taylor series to confirm this at small times, $t \ll 1$, and then the semigroup property brings confirmation for all time $t>0$.
When $\alpha \ne 1$, \eqref{eq:mittag:leffler:negative:argument:integral:1:matrix:argument} is the Mittag-Leffler function of the matrix.
Notice that $(1,1,\ldots,1) \mathbb{A} = (0,0,\ldots,0)$ and $(1,1,\ldots,1) \mathbb{I} = (1,1,\ldots,1)$, so
multiplication of the Taylor series on the left by a row vector of ones quickly shows that columns of $E_{\alpha}(\mathbb{A}t^{\alpha})$ always sum to one.
 If we could show that all entries of the Mittag-Leffler function of a graph Laplacian matrix  are nonnegative (and they are!) then we would know that $E_{\alpha}(\mathbb{A}t^{\alpha})$ is also a stochastic matrix.
Unlike the exponential case though, we do not have the semigroup property.
Nevertheless the Taylor series still confirms nonnegativity for \textit{small times}, because $E_{\alpha}(\mathbb{A}t^{\alpha}) \approx \mathbb{I} + t^{\alpha} \mathbb{A} + \mathcal{O}(t^{2 \alpha})$, off-diagonals are nonnegative, and for small $t$ the identity matrix more than compensates for the negative diagonal of $\mathbb{A}$.


\textit{A $2 \times 2$   matrix example.}
For $a,b>0$, the diagonalization
 \begin{equation}
\left(
\begin{tabular}{rr}
$-a$ & $b$ \\
$a$  & $- b$
\end{tabular}
\right)
=
\left(
\begin{tabular}{rr}
$b$ & $1$ \\
$a$  & $-1$
\end{tabular}
\right)
\left(
\begin{tabular}{rr}
$0$ & $0$ \\
$0$  & $a+b$
\end{tabular}
\right)
\,
\left(
\begin{tabular}{rr}
$-1$ & $-1$ \\
$-a$  & $b$
\end{tabular}
\right)
\frac{1}{a+b} \nonumber
\end{equation}
shows  $-(-\mathbb{A})^{\frac{1}{\alpha}}=  \mathbb{A} (a+b)^{\frac{1}{\alpha} -1} $ is still a graph Laplacian
\[
-(-\mathbb{A})^{\frac{1}{\alpha}} =
\left(
\begin{tabular}{rr}
$-a$ & $b$ \\
$a$  & $- b$
\end{tabular}
\right)^{\frac{1}{\alpha}}
=
\left(
\begin{tabular}{rr}
$-a$ & $b$ \\
$a$  & $- b$
\end{tabular}
\right)
(a+b)^{\frac{1}{\alpha} -1} .
\]
Analogous to the  scalar fractional Euler generalisation  \eqref{eq:weighted:sum:euler:limits:negative}, to accommodate a matrix  \eqref{eq:weighted:sum:euler:limits:negative} could become a discretization of \eqref{eq:mittag:leffler:negative:argument:integral:1:matrix:argument}:
\begin{equation}
\sum_k w_k   \left( \mathbb{I} + \frac{-s_k(-\mathbb{A})^{\frac{1}{\alpha}}t}{n} \right)^n .
\label{eq:weighted:sum:euler:limits:negative:matrix:argument}
\end{equation}
Here $w_k = w_{-,1}(s_k) \ge 0$ is a discretization of \eqref{eq:mittag:leffler:density:w:minus:1}.
For all sufficiently large $n$ the entries of $ ( \mathbb{I} -s_k(-\mathbb{A})^{\frac{1}{\alpha}}t/n )$ are nonnegative.  
Thus we see via \eqref{eq:weighted:sum:euler:limits:negative:matrix:argument} or \eqref{eq:mittag:leffler:negative:argument:integral:1:matrix:argument}  that the Mittag-Leffler function is a stochastic matrix.
Similar to the numerical method of \textit{uniformization}, the matrix $( \mathbb{I} -s_k(-\mathbb{A})^{\frac{1}{\alpha}}t/n )$ can be interpreted as a first order approximation of a Markov process in continuous time by a  Markov process in discrete time.
We can interpret   $E_{\alpha}(\mathbb{A} t^{\alpha})$ (which may itself be a \textit{semi-Markov process} \cite{SemiMarkovGraph2011}) as a weighted sum of related Markov processes.
Denote \eqref{eq:mittag:leffler:density:w:plus:big:W} with $\lambda=1$ by $W_{+,1}(s)$ and rearrange  \eqref{eq:mittag:leffler:positive:argument:integral} to
$\mathbb{P} =  \frac{\alpha}{1-\alpha} ( \frac{\exp((-\mathbb{A})^{\frac{1}{\alpha}} t)}{\alpha} - E_{\alpha}(- \mathbb{A} t^{\alpha})  ) = \int_0^{\infty}  W_{\bm{+,1}}(s)  \exp(-s (-\mathbb{A})^\frac{1}{\alpha} t) \ud s.$
If $-(-\mathbb{A})^\frac{1}{\alpha}$ is also a graph Laplacian then  similarly $\mathbb{P}$ must also be a stochastic matrix.



\subsection{A Mittag-Leffler function of a graph Laplacian is a stochastic matrix}
Instead of the explicit, forward difference \eqref{eq:Euler:finite:difference} that led to our first candidate, we could replace \eqref{eq:Euler:finite:difference} by an implicit, \textit{backward difference}
of the same continuous equation, $\ud y/ \ud t =  y$.
(In fact the backward differences of the Grunw\"ald-Letnikov approach in Section~\ref{sec:Grunwald:Letnikov} are one avenue to the fractional resolvent that we now describe \cite{BajlekovaThesis2001}.)
That  leads to the discrete approximation
$
(y_{j}-y_{j-1})/h =  y_{j}.
$
Set $h= t/n$.
Recursively,
$
 y_n =   \left(1 -   t/n  \right)^{-n}.
$
This is another Euler formula, similar to that in Table~\ref{tab:exp:mittag:leffler}.
Both converge to the same exponential limit,
\begin{equation}
\lim_{n \rightarrow \infty} \left(1 +  \frac{t}{n}  \right)^{n} \;=\; \lim_{n \rightarrow \infty} \bm{\left(1 \bm{-}  \frac{t}{n}  \right)^{\bm{-}n}} \;=\; \exp(t).
\label{eq:euler:limit:forward:backward}
\end{equation}
However, the latter representation is more suggestive of the important connection to the \textit{resolvent}  $(s\mathbb{I} - \mathbb{A})^{-1}$.
This resolvent matrix at $s$ is defined whenever $s$ is not an eigenvalue, and is involved in solutions of  $\ud y/ \ud t =  \mathbb{A}y$, via an Euler limit
\begin{equation}
\left( \mathbb{I} - \frac{t}{n}\mathbb{A} \right)^{-n} \longrightarrow \exp(\mathbb{A}t).
\label{eq:euler:limit:matrix:exponential}
\end{equation}

Central to the Hille-Yosida and related theorems for semigroups is the relationship between the resolvent and the  exponential of the infinitesimal generator $\mathbb{A}$.
Namely, the resolvent is the Laplace transform of the exponential: $  (s \mathbb{I} - \mathbb{A})^{-1} = \int_0^{\infty} \exp(-st) \exp(\mathbb{A}t) \ud s.$
The Mittag-Leffler functions do not have the semigroup property but the essence of the relationship can be generalised to express the Laplace transform of the Mittag-Leffler function, $E_{\alpha} (\mathbb{A} t^{\alpha})$, in terms of the resolvent:
\begin{equation}
s^{\alpha -1} (s^{\alpha} \mathbb{I} - \mathbb{A})^{-1} = \mathcal{L} \left\{ E_{\alpha} (\mathbb{A} t^{\alpha})  \right\}  =  \int_0^{\infty} e^{-st} E_{\alpha} (\mathbb{A} t^{\alpha}) \ud s.
\label{eq:resolvent:is:laplace:transform:of:mittag:leffler}
\end{equation}
The scalar version of \eqref{eq:resolvent:is:laplace:transform:of:mittag:leffler} appears in \eqref{eq:laplace:transform:mittag:leffler:scalar} and very similar steps lead to \eqref{eq:resolvent:is:laplace:transform:of:mittag:leffler}.
Thus we can now take the inverse Laplace transform of \eqref{eq:resolvent:is:laplace:transform:of:mittag:leffler} to express the Mittag-Leffler function in terms of the resolvent.
This could lead to what amounts to the same representation as in \eqref{eq:Mittag:Leffler:contour:integral:representation:lambda:alpha}:
$
 E_{\alpha}(\mathbb{A} t^{\alpha}) = \frac{1}{2 \pi i}  \int_{\mathcal{C}} e^{s} s^{\alpha-1} (s^{\alpha}\mathbb{I} -\mathbb{A} t^{\alpha})^{-1} \ud s.
$
But a different representation appears by instead using the \textit{Post--Widder inversion formula}:
\[
f(t) = \lim_{n \rightarrow \infty} \frac{(-1)^n}{n!} \left( \frac{n}{t} \right)^{n+1} \left( \frac{\ud^n}{\ud s^n} \hat{f} \right) \left(\frac{n}{t} \right).
\]
That inversion formula comes from the usual rule for Laplace transforms that the $nth$ derivative  in the $s$-domain, $\hat{f}^{(n)}(s)$, is paired with  $(-1)^n t^n f(t)$ in the time domain, and by noticing that the corresponding Laplace integral transforms $\int_0^\infty \rho_n(s) f(s) \ud s \rightarrow f(t)$, because the $\rho_n$ tend to the Dirac delta distribution.
This leads to Bajlekova's representation \cite[Proposition 2.11]{BajlekovaThesis2001} of  $E_{\alpha}(\mathbb{A} t^{\alpha})$:
\begin{equation}
\lim_{n \rightarrow \infty} \frac{1}{n!} \sum_{k=1}^{n+1} b_{k,n+1} \left( \mathbb{I} - (t/n)^{\alpha}\mathbb{A} \right)^{-k}.
\label{eq:fractional:resolvent}
\end{equation}
The $b_k$ are the positive constants 
in the $n$th derivative     $\ud^n / \ud s^n\, ( s^{\alpha-1} \left(s^{\alpha} \mathbb{I} - \mathbb{A} \right)^{-1} ) =
 (-1)^n s^{-n-1} \sum_{k=1}^{n+1} b_{k,n+1} \left( s^{\alpha} \left( s^{\alpha}\mathbb{I} - \mathbb{A} \right) \right)^{-k}
$
that Post--Widder requires.

Representation \eqref{eq:fractional:resolvent} is a generalisation of \eqref{eq:euler:limit:matrix:exponential} and in that sense it is yet another fractional generalisation of the Euler formula.
Although \eqref{eq:fractional:resolvent} would not usually be a good numerical scheme, it can be usefully  applied to affirmatively answer our earlier question concerning nonnegativity of the Mittag-Leffler function of a graph Laplacian, $E_{\alpha}(\mathbb{A} t^{\alpha})$.
First note the pattern of $\pm$ signs of the graph Laplacian $\mathbb{A}$ implies the following simple pattern of signs in $\left( \mathbb{I} - (t/n)^{\alpha}\mathbb{A} \right)$: positive entries on the main diagonal, and negative entries off the main diagonal (although possibly zero entries are allowed).
That pattern is displayed here on the left of \eqref{eq:sign:pattern:Euler} for the $3 \times 3$ case:
\begin{equation}
\left(
\begin{array}{ccc}
\bm{+} & - & - \\
- & \bm{+} & - \\
- & - & \bm{+}
\end{array}
\right)^{-1}
= \;\;
\left(
\begin{array}{ccc}
+ & + & + \\
+ & + & + \\
+ & + & +
\end{array}
\right).
\label{eq:sign:pattern:Euler}
\end{equation}
If a matrix with this pattern of signs has sufficiently large entries on the main diagonal, then all entries of the inverse of that matrix are positive.
That `inverse positive' property displayed schematically in \eqref{eq:sign:pattern:Euler} is a fact of linear algebra \cite{Str09}, related to the class of $M-$matrices \cite{Berman:1994aa}, and it has a  generalisation to operators.
Applied to our examples (for all sufficiently large $n$ when the diagonal entries are relatively large enough) the pattern of signs \eqref{eq:sign:pattern:Euler} implies that \textit{all} entries of  $\left( \mathbb{I} - (t/n)^{\alpha}\mathbb{A} \right)^{-1}$ are nonnegative.
Powers of a nonnegative matrix are of course nonnegative so  $\left( \mathbb{I} - (t/n)^{\alpha}\mathbb{A} \right)^{-k}$ is likewise nonnegative.
Representation \eqref{eq:fractional:resolvent} now shows that all entries of $E_{\alpha}(\mathbb{A} t^{\alpha})$ are nonnegative.
Having already established unit column sum, we have shown:  \textit{the Mittag-Leffler function of a graph Laplacian, $E_{\alpha}(\mathbb{A} t^{\alpha})$, is indeed a stochastic matrix.}

\subsection{Random walks on a graph Laplacian}
A continuous time random walk (CTRW) \cite{KlafterSokolovBook2011} on the set of  states $\{1,2, \ldots,n\}$ is associated with the  $n \times n$ graph Laplacian matrix $\mathbb{A}$ as follows.
Having just arrived at a state $j$, the walker waits in that state for an amount of time that is a random variable with a state-dependent probability density $\psi(j,t)$.
We have  $\psi(j,t)>0$ and $1 = \int_0^{\infty} \psi(j,t) \ud t $ because $\psi$ is a density.
The walker then immediately moves to a different state $i$  with probability
\begin{equation}
\lambda(i,j) = \frac{a_{ij}}{\left| a_{jj} \right|}  \qquad (i \ne j).
\label{eq:lambda:ij:a:ij}
\end{equation}
We also define $\lambda(j,j)= \lambda_{jj} = \left| a_{jj} \right|$ (unlike the convention $\lambda_{jj}=0$ of some authors).
Notice that the properties of the matrix ensure that $\lambda(i,j)>0$ and $1 = \sum_{i, i \ne j} \lambda(i,j)$.

\textit{The memory function.}
This CTRW is a Markov process if and only if the waiting time in all states is an exponential random variable, so that $\psi(j,t) = a_{jj} \exp(-a_{jj} t)$.
 In that special case the master equation that governs the evolution of the probability $p(j,t)$ of being in state $j$ at time $t$ is $\frac{\ud }{\ud t} \bm{p}(t)  =  \mathbb{A} \bm{p}(t)$ with solution $\bm{p}(t)=\exp(\mathbb{A}t) \bm{p}(0)$, where $\bm{p}(t)$ is a vector with entries $p(j,t)$.
Next, we examine the master equation that governs the probability associated with this CTRW in the more general case that does not make assumptions about the form of the waiting time density --- it need not be exponential for example, in which case the process is not Markov and must therefore exhibit some form of memory.
That notion of memory turns out to be made mathematically precise by the ratio of the Laplace transform of the survival function to that of the waiting time density:
\begin{equation}
\hat{K}(s) = \frac{\hat{\psi}(s)}{\hat{\phi}(s)} =  \frac{s\hat{\phi}(s)-1}{\hat{\phi}(s)}.
\label{eq:memory:function}
\end{equation}
Mainardi \textit{et al.} termed \eqref{eq:memory:function} the \textit{memory function} of the process (although they work with the reciprocal of our notation) \cite{MainardiRabertoGorenfloScalas2000}, and it characterises the generalised master equation
\eqref{eq:master:equation:general} that we describe next.

\subsection{A generalised master equation for waiting times that are not exponential}

It is difficult to generalise to arbitrary waiting times by working with the probability density description alone.
Others   have overcome this by a finer description of the process in terms of a flux, and we follow a very similar derivation here.
This use of a flux after $n$ steps is analogous to the way subordinators and random time change representations are used in stochastic process theory to connect absolute time to the random number of steps that occur.
It is also analogous to derivations of the classical heat equation that do not simply work with temperature alone, and require a notion of heat flux \cite{Fedotov:2010aa,SokolovSchmidtSagues2006,AngstmannDonnellyHenry2013,ChechkinGorenfloSokolov2005,Kur80, PrussBook2012}, \cite[Ch. 5]{KlafterSokolovBook2011}.

Define
$q_n(j,t) \equiv $  \textit{probability density to arrive at site  $j$  at time  $t$,  after exactly     $n$  steps.}
This nonnegative flux measures the `rate of flow of probability \textit{in}' to the state.
Initially $(t=0)$ assume the particle is in state $j_0$, and that the particle only just arrived, so \;
$
q_0(j,t) \equiv \delta_{j,j_0} \delta(t)
$
where $\delta$ denotes Dirac's delta distribution.
Sum over all steps to get a total flux \;
$
q(j,t) \equiv \sum_{n=0}^{\infty}   q_n(j,t) .
$
It is convenient to have notation for the total flux minus the initial flux
\begin{equation}
q^{+}(j,t) \equiv \sum_{n\bm{=1}}^{\infty}   q_n(j,t) \; = q(j,t) - q_0(j,t) = q(j,t) - \delta_{j,j_0} \delta(t).
\label{eq:q:plus}
\end{equation}
Recursively, for $n \ge 1$,
\begin{equation}
q_{n+1}(j,t) = \sum_{i, i \ne j} \int_0^t  q_n(i,u) \Psi(j,t,i,u) \ud u
\label{eq:q:n:recursive}
\end{equation}
where
$\Psi(j,t,i,u) \equiv $
\textit{probability density to arrive at site $j$ at time  $t$ in one step, given  it arrived at site $i$ at time $u$.}
Put \eqref{eq:q:n:recursive} in \eqref{eq:q:plus} and swap order of summation
$q^{+}(j,t)=\sum_{n=1}^{\infty}  ( \sum_{i \ne j} \int_0^t  q_{n-1}(i,u) \Psi(j,t,i,u) \ud u )   =\sum_{i \ne j} \int_0^t  (\sum_{n=1}^{\infty}q_{n-1}(i,u) )  \Psi(j,t,i,u) \ud u$.
If we substitute
$ \sum_{n=1}^{\infty}q_{n-1}(i,u)   = \sum_{n=0}^{\infty}q_{n}(i,u)  =     q(i,u) $
and assume the functional form
$
\Psi(j,t,i,u) = \lambda(j,i) \psi(i,t-u)
$
where
$\psi(i,t) \equiv $
\textit{probability density to leave site $i$ at time  $t$, given it arrived at site $i$ at $t=0$}
 is the \textit{waiting time}, then
 \begin{equation}
  q^{+}(j,t)=  \sum_{i \ne j} \lambda(j,i) \int_0^t  q(i,u)  \psi(i,t-u) \ud u.
  \label{eq:q:plus:2}
\end{equation}

The probability density $p(j,t)$ to be at state $j$ at time $t$, is related to the flux by
\begin{equation}
p(j,t) = \int_0^t \phi(j,t-u) q(j,u) \ud u
\label{eq:p:convolution:identity}
\end{equation}
where the \textit{survival time}
$\phi(j,t) \equiv$  \textit{the probability to remain at site $j$ for all times in $(0,t)$, having just arrived at  $t=0$.}
Using \eqref{eq:q:plus} gives
$
p(j,t) =  \delta_{j,j_0}  \phi(j,t) + \int_{0^+}^t \phi(j,t-u) q^{+}(j,u) \ud u
$
where the integral with lower limit $0^{+}$ does not include $0$.
Differentiating,
\begin{eqnarray}
\frac{\ud }{\ud t}  p(j,t) &=& \delta_{j,j_0}   \frac{\ud }{\ud t}\phi(j,t)  +\; \frac{\ud }{\ud t} \int_{0^+}^t \phi(j,t-u) q^{+}(j,u) \ud u \nonumber \\
               &=& \delta_{j,j_0}   (-\psi(j,t))  +   \phi(j,t-t) q^{+}(j,t) + \int_{0^+}^t (-\psi(j,t-u)) q^{+}(j,u) \ud u \nonumber \\
                   & =&  - \int_{0^+}^t \psi(j,t-u) q(j,u) \ud u  + \sum_{i \ne j} \lambda(j,i) \int_0^t  q(i,u)  \psi(i,t-u) \ud u.
                  \label{eq:p:derivative}
\end{eqnarray}
Note $\phi(j,t-t)=\phi(j,0)=1$.
We also used \eqref{eq:q:plus}, \eqref{eq:q:plus:2} and the multivariate chain rule to differentiate the integral.

We now have two equations involving the probability density: \eqref{eq:p:convolution:identity} and \eqref{eq:p:derivative}.
Take the Laplace transform of both.
The first equation, \eqref{eq:p:convolution:identity}, is a convolution so the transform is a product of transforms: \,
$
\hat{p}(j,s) = \hat{\phi}(j,s)\hat{q}(j,s).
$
Multiply by $\hat{\psi}(j,s)$ and divide by $\hat{\phi}(j,s)$ to get
\begin{equation}
  \hat{q}(j,s) \hat{\psi}(j,s) =   \frac{ \hat{\psi}(j,s)}{\hat{\phi}(j,s)} \hat{p}(j,s) = \hat{K}(j,s) \hat{p}(j,s)
\label{eq:replacement}
\end{equation}
where the state-dependent version of the memory function \eqref{eq:memory:function} now appears
\begin{equation}
\hat{K}(j,s) \equiv \frac{ \hat{\psi}(j,s)}{\hat{\phi}(j,s)} .
\label{eq:K:Laplace:transform}
\end{equation}
The second equation \eqref{eq:p:derivative} involves a derivative on the left and convolutions on the right so the Laplace transform is
\[
s \hat{p}(j,s) - p(j,0) =   -  \hat{\psi}(j,s) \hat{q}(j,s)
              + \sum_{i \ne j} \lambda(j,i)   \hat{\psi}(i,s) \hat{q}(i,s).
\]
Use \eqref{eq:replacement} and \eqref{eq:K:Laplace:transform} to replace $ \hat{\psi}(i,s) \hat{q}(i,s)$ by $\hat{K}(i,s) \hat{p}(i,s)$:
\begin{equation}
s \hat{p}(j,s) - p(j,0) =   - \hat{K}(j,s) \hat{p}(j,s)
              + \sum_{i \ne j} \lambda(j,i)  \hat{K}(i,s) \hat{p}(i,s).
\label{eq:Laplace:transform:K:me}
\end{equation}
Take inverse Laplace transforms to finally arrive at the desired master equation
\begin{eqnarray}
 \frac{\ud }{\ud t}  p(j,t)  &=&   -  \int_0^{t} K(j,t-u) p(j,u)  \ud u + \sum_{i \ne j} \lambda(j,i) \int_0^{t} K(i,t-u) p(i,u)  \ud u.
\label{eq:master:equation:general}
\end{eqnarray}
This master equation does not assume exponential waiting times.
The waiting times may have different functional forms in different states.
In the special case of exponential waiting times, $\hat{K}(j,s)  = \lambda_{jj}$ in \eqref{eq:Laplace:transform:K:me} so the $K(j,t)$ appearing in the convolutions in \eqref{eq:master:equation:general} are Dirac delta distributions and \eqref{eq:master:equation:general} collapses to the usual master equation $\frac{\ud }{\ud t} \bm{p}(t)  =  \mathbb{A} \bm{p}(t)$.

\subsection{Mittag-Leffler waiting times have Mittag-Leffler matrix functions for master equations}

 Now specialise the general master equation \eqref{eq:master:equation:general} to  Mittag-Leffler waiting times:  state $j$ has survival function $\phi(j,t) = E_{\alpha}(-\lambda_{jj}t^{\alpha})$ \cite{HilferAnton95}.
The following steps are very similar to \cite{MainardiRabertoGorenfloScalas2000}, although here we  work with a matrix.
By \eqref{eq:laplace:transform:mittag:lefflerminus:lambda},
\[
\hat{\phi}(j,s) =  \mathcal{L} \left\{ \mathcal{L} \left\{ w_{-}(s) \right\} \right\} = \mathcal{L} \left\{ E_{\alpha}(-\lambda_{jj}t^{\alpha}) \right\} = \frac{s^{\alpha-1}}{\lambda_{jj}+s^{\alpha}}.
\]
The transform of its derivative $\psi = - \frac{d}{dt} \phi$ is $s\hat{\phi}(s)-1$ so
$
\hat{\psi}(j, s) = \frac{\lambda_{jj}}{\lambda_{jj}+s^{\alpha}} .
$
The memory function   \eqref{eq:K:Laplace:transform} is thus:
\[
\hat{K}(j,s)  = \lambda_{jj} s^{1- \alpha}  .
\]
Substitute into  \eqref{eq:Laplace:transform:K:me}:
$
s \hat{p}(j,s) - p(j,0) =   -\lambda_{jj} s^{1- \alpha}  \hat{p}(j,s)
              + \sum_{i \ne j} \lambda(j,i)  \lambda_{ii} s^{1- \alpha} \hat{p}(i,s).
$
Then divide both sides by $s^{1- \alpha}$, so
\begin{equation}
s^{\alpha} \hat{p}(j,s) - s^{\alpha-1}p(j,0) =   -\lambda_{jj}   \hat{p}(j,s)
              + \sum_{i \ne j} \lambda(j,i)  \lambda_{ii}  \hat{p}(i,s).
 \label{eq:mittag:leffler:master:equation:laplace:transform}
\end{equation}
The right side is precisely the matrix-vector product $\mathbb{A} \bm{\hat{p}}$ (recalling \eqref{eq:lambda:ij:a:ij}).
We recognise the left as the \textit{Laplace transform of a Caputo fractional derivative}:
\[
\mathcal{L} \left\{  \frac{\ud^{\alpha} }{\ud t^{\alpha}} f(t)    \right\} = s^{\alpha} \hat{f}(s) -  s^{\alpha-1} f(0^{+}).
\]
Thus, after an inverse transform, the master equation is
$
\frac{\ud^{\alpha} }{\ud t^{\alpha}} p(t)  =  \mathbb{A} p(t) .
$
We may also write this as the matrix-vector version of the scalar model fractional equation \eqref{eq:fractional:ode} from the Introduction (although \eqref{eq:mittag:leffler:master:equation} corresponds to the fractional decay case)
\begin{equation}
D_t^{\alpha} \bm{p} = \mathbb{A} \bm{p}  \qquad \textrm{\;\; with solution} \quad \bm{p}(t) = E_{\alpha}(\mathbb{A} t^{\alpha}) \bm{p}(0).
\label{eq:mittag:leffler:master:equation}
\end{equation}
Conclusion: \textit{Mittag-Leffler waiting times go together with a Mittag-Leffler function of the (graph Laplacian) matrix.}
This is also a second proof of our earlier observation near \eqref{eq:sign:pattern:Euler} that a Mittag-Leffler function of a graph Laplacian is a stochastic matrix.

\subsection{A remarkable class of nonnegative matrix functions}
While  \eqref{eq:master:equation:general} boasts generality,
 it does not make apparent the simple structure of important special classes of solutions.
Let us now elucidate one such special class, in terms of matrix functions.
The algebra that permits the simplified form of \eqref{eq:mittag:leffler:master:equation:laplace:transform} and thus \eqref{eq:mittag:leffler:master:equation}, is the old idea of separation of variables: $\hat{K}(j,s)  = \hat{k}(s) g(j) $.
The generalised master equation then reads
\[
\frac{1}{\hat{k}(s)} \left( s \hat{p}(j,s) - p(j,0) \right) =   - g(j)   \hat{p}(j,s)
              + \sum_{i \ne j} \lambda(j,i)  g(i) \hat{p}(i,s).
\]
The right is $\mathbb{A} \bm{\hat{p}}$ where the main diagonal entries are $a_{jj} = - g(j)$.
In the time domain, the left is
\begin{equation}
\int_0^t \frac{1}{k(t-u)} \frac{\ud p}{\ud u}   \ud u.
\label{eq:convolution:memory:function}
\end{equation}
The Caputo derivative \eqref{eq:Caputo:derivative:definition} corresponds to a well-understood example
\begin{equation}
k(t) \equiv t^{\alpha} \Gamma(1-\alpha),
\label{eq:Caputo:memory:function}
\end{equation}
and the memory terms appearing in our first candidate for the fractional Euler formula \eqref{eq:fractional:Euler:limit} can now be understood as coming from a discretization of the corresponding convolution \eqref{eq:convolution:memory:function} with that particular form of $k$.

Solutions are matrix functions, computed as the inverse Laplace transform of
\begin{equation}
\frac{1}{\hat{k}(s)} \left( \frac{s}{\hat{k}(s)}  \mathbb{I} - \mathbb{A} \right)^{-1}.
\label{eq:Laplace:transform:solution:separable:case}
\end{equation}
Numerical inversion of \eqref{eq:Laplace:transform:solution:separable:case} is well suited to Cauchy contour integral methods, although further research is required to address issues associated with \textit{pseudospectra} of graph Laplacians
\cite{ShevCauchyIntegralMasterEqnPseudoSpectraCTAC2015,TreEmb05}, and to tailor those methods to various forms of the memory function $\hat{k}$.
As probabilities, these solutions are nonnegative, so exploring separable forms of the memory function gives rise to a large class of nonnegative matrix functions.
Indeed, these functions output stochastic matrices.
Confirming unit column sum might again come by multiplying a Taylor series in powers of $\mathbb{A}$ on the left by a row  vector of ones.
That the first term in such a series should be simply $\mathbb{I}$ might be established by applying the usual Laplace transform rule that $f(0^+) = \lim_{s \rightarrow \infty} s \hat{f}(s)$ to \eqref{eq:Laplace:transform:solution:separable:case}.
More interestingly, the very special nonnegativity property is explained as variations of the algebraic structure elucidated in  \eqref{eq:fractional:resolvent} and \eqref{eq:sign:pattern:Euler}.
Inversion of the Laplace transform \eqref{eq:Laplace:transform:solution:separable:case} could express the solution in terms of a resolvent, and then
the key observation displayed in \eqref{eq:sign:pattern:Euler} is that our graph Laplacian is an example of a \textit{resolvent positive operator} \cite{AltenbergPNASResolventPositive2012,LaplaceTransformCauchyProblemBook2011}: for large positive $s$, \textit{all entries of the resolvent matrix $(s\mathbb{I} - \mathbb{A})^{-1}$ are nonnegative.}
The resolvent shares the completely monotone properties of section \ref{sec:completely:monotone}, and the usual Euler formula in terms of the resolvent \eqref{eq:euler:limit:matrix:exponential} quickly shows solutions of the corresponding Cauchy problem preserve positivity.
Further exploration of fractional Euler formulas (as variations of \eqref{eq:weighted:sum:euler:limits:negative:matrix:argument} or \eqref{eq:fractional:resolvent}), and further exploration of algebraic structures of master equations including Toeplitz or Hankel structures and other splittings \cite{StrMac14,SpeGre13}, therefore seems promising.

\begin{figure}[ht]
\centering
\begin{tabular}{cc}
\includegraphics[scale=0.35]{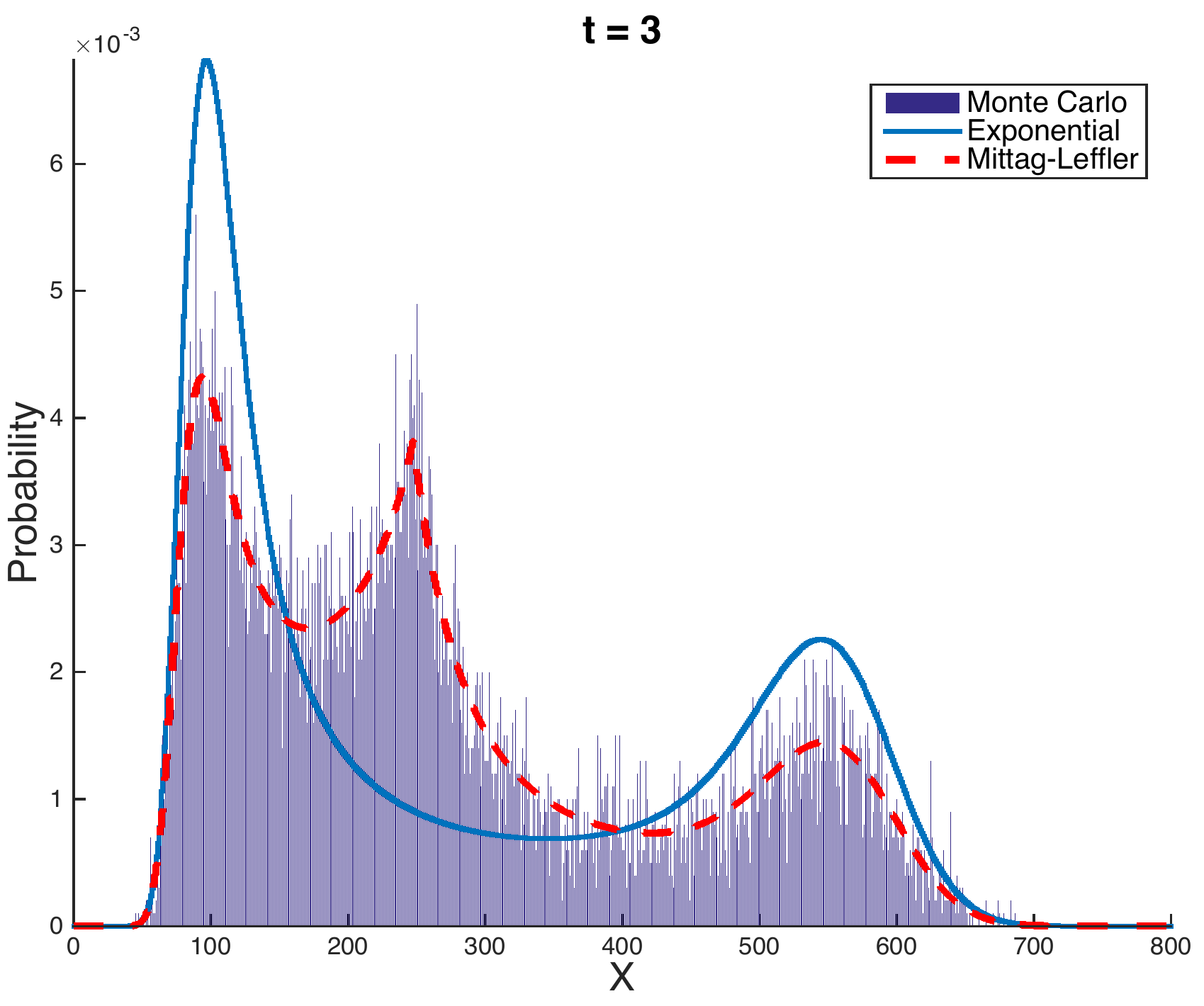} &
\includegraphics[scale=0.3]{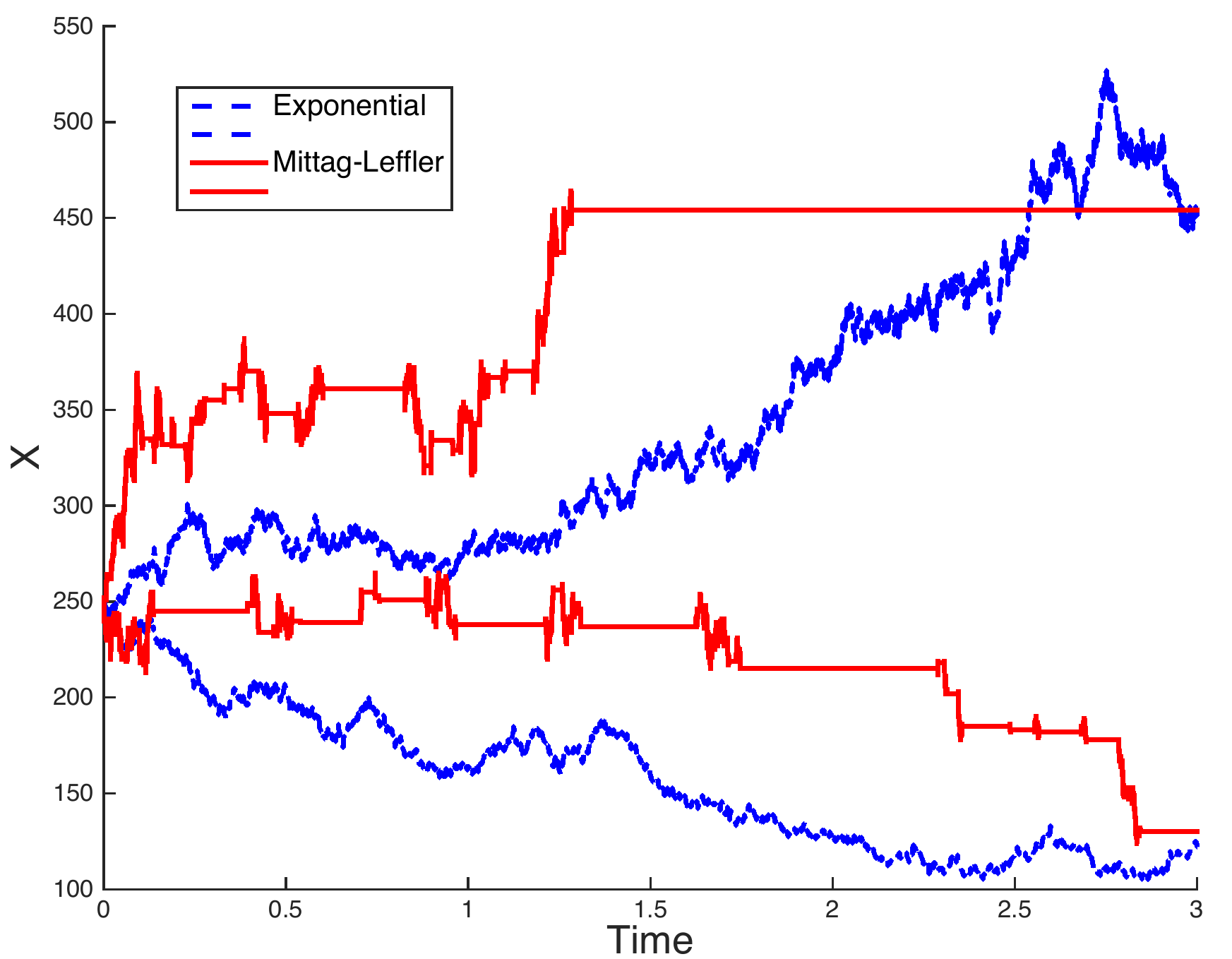}
\end{tabular}
\caption{\textit{Left:} The Schl\"ogl reactions with Mittag-Leffler waiting times have a Mittag-Leffler solution \eqref{eq:mittag:leffler:master:equation} coming from the generalised master equation \eqref{eq:master:equation:general}.
This is distinctly different to the solution of the usual master equation with exponential waiting times.
\textit{Right:} Gillespie-like stochastic simulations of the Schl\"ogl reactions, with the usual exponential waiting times, and also with Mittag-Leffler waiting times \eqref{eq:inverse:CDF:w:minus:mixture:exponential:Gillespie}.}
\label{fig:ML}
\end{figure}


 \section{\label{sec:Schlogl} Application to the Schl\"ogl reactions}

The Sch\"ogl model consists of four reactions
\begin{center}
\begin{tabular}{ccccc}
1. & $B_1 + 2X  \rightarrow 3X$ & & & $k_1 = 3 \times 10^{-7}$ \\
2. & $B_1 + 2X  \leftarrow 3X$ && & $k_2 = 1 \times 10^{-4}$ \\
3. & $B_2 \rightarrow X$ & & & $k_3 = 1 \times 10^{-3}$\\
4. & $B_2 \leftarrow X$   & & & $k_4 = 3.5$
\end{tabular}
\end{center}
Here $B_1=1 \times 10^{5}$ and $B_2=2 \times 10^{5}$ are constants that model buffered species.
We choose the initial condition $X(0) = 247$, which in the deterministic version of the model lies on the separatrix between the two basins of attraction of the lower steady state (at about $85$) and of the higher steady state (at about $565$).
The stochastic model exhibits a bimodal distribution (Figure~\ref{fig:ML}).
Fedatov also uses this model as an example, in a physically comprehensive discussion \cite{Fedotov:2010aa}.

Such reactions can be modelled in the framework of the chemical master equation, which is a memoryless Markov process with exponential waiting times \cite{Gil92,PetzoldCloud}.
Waiting times different from the usual exponential choice could model some form of anomalous diffusion, although there are many issues associated with a careful physical interpretation \cite{Fedotov:2010aa, HellanderLotstedt2016}.
Operationally, experimentation with different waiting times comes by simply changing the simulation of the waiting time in the  usual Gillespie stochastic simulation algorithm, while keeping all other steps in the algorithm the same \cite{Gil92}.
 As an example, Figure \ref{fig:ML} incorporates Mittag-Leffler waiting times  ($\alpha = 0.7$) in simulations of the Schl\"ogl reactions, so the generalised master equation is  \eqref{eq:mittag:leffler:master:equation}.
The waiting time between reactions is simulated via \eqref{eq:sampling:mittag:leffler}:
\begin{equation}
\tau \; \sim \;\;\; - \left(\frac{1}{a}\right)^\frac{1}{\alpha}  \left( \frac{\sin(\pi \alpha)}{\tan(\pi \alpha (1-u_1))} - \cos(\pi \alpha) \right)^\frac{1}{\alpha}  \log (u_2)
\label{eq:inverse:CDF:w:minus:mixture:exponential:Gillespie}
\end{equation}
where $u_1,u_2$ are drawn from independent uniform random variables on $(0,1)$, and where $a$ is the usual sum of the propensities of the reactions in the current state.

For this model there are some numerical issues associated with truncation to a finite state space \cite{MacBur08} and with stationary distributions \cite{Lythe2009}, although numerical experiments indicate good accuracy here.
 The Mittag-Leffler waiting times manifest themselves in the very long pauses between consecutive reactions (Figure~\ref{fig:ML}).
  With the usual exponential waiting times the solution of the associated master equation is only bimodal, but the Mittag-Leffler solution exhibits a third mode at $\approx 247$ that persists for a very long time (Figure~\ref{fig:ML}).

The generalised master equation \eqref{eq:master:equation:general} invites further experimentation with different waiting times leading to processes with memory.
One form of memory could come from the $W_{+}$ mixture in \eqref{eq:W:plus:exponential:mixture}.
A fast Monte Carlo method to sample from that $W_+$ mixture is provided in \eqref{eq:sampling:W:plus:mixture:of:exponentials}.
The generalised master equation would be \eqref{eq:master:equation:general} with a particular memory function  $\hat{K}_{W_{+}}(s) = \hat{\phi}_{W_{+}}(s)/(s  \hat{\phi}_{W_{+}} (s)-1)$, with transform $ \hat{\phi}_{W_{+}}$ available in \eqref{eq:laplace:transform:W:plus}.

\section*{Discussion}
Our fractional usury was motivated in part by a quest for the `right' notion of \textit{fractional compound interest}, a subject that deserves more attention.
With this in mind various candidates for a  fractional generalisation of the closely related Euler formula  have been discussed.
Three in particular were based on: generalising the discrete construction that leads to the usual Euler limit \eqref{eq:fractional:Euler:limit}, a Gr\"unwald-Letnikov approach \eqref{eq:GL:fractional:Euler:limit},  or a Cauchy integral \eqref{eq:weighted:sum:euler:limits:negative}-\eqref{eq:weighted:sum:euler:limits:positive}.
That all candidates be discrete and converge to the continuous  Mittag-Leffler function has served as a guiding principle.
Both a Cauchy formulation and a formulation in terms of the resolvent are attractive because they have the satisfying form of being a weighted sum of regular Euler formulas.
Together with the observation that the graph Laplacian is a resolvent positive operator, with pattern of signs as in \eqref{eq:sign:pattern:Euler}, one application of these fractional Euler formulas is to show that the Mittag-Leffler function of a graph Laplacian is a stochastic matrix.
Generalisations of other exponential properties to the Mittag-Leffler setting are destined.
For instance, the Golden-Thompson inequality \cite{ForresterGoldenThompson2014}  states that, for Hermitian matrices $A$ and $B$,
$
\textrm{trace} \left( \exp(A+B) \right) \le \textrm{trace} \left( \exp(A) \exp(B) \right).
$
The inequality  fails when the exponential function is replaced by the Mittag-Leffler function, as a scalar example quickly shows $E_{1/2}(1+1) = \exp((1+1)^2) \textrm{erfc}(-(1+1)) \approx 109 \not \leq  5 \times 5 \approx E_{1/2}(1) \times E_{1/2}(1)$, but less na\"ive generalisations might be possible.
Finally, continuing this theme, we have elucidated connections to master equations that generalise the usual exponential waiting times to Mittag-Leffler waiting times, exemplified by applications to modeling and simulation of chemical reactions.


\bibliographystyle{siamplain}
\bibliography{shevrefs}
\end{document}